\documentclass[aos,MSNbibl,citesort,seceqn,dvips]{arximspdf}
\usepackage{graphicx}
\doi{10.1214/10-AOS868}
\volume{39}
\issue{3}
\pubyear{2011}
\firstpage{1427}
\lastpage{1470}

\makeatletter

\newtheorem{theorem}{Theorem}[section]
\newtheorem{lemma}{Lemma}[section]
\newtheorem{cor}{Corollary}[section]
\newtheorem{prop}{Proposition}[section]

\newproclaim{rem}{Remark}[section]

\newproclaim{assumptionP}{Assumption}
\newproclaim{assumptionA}{Assumption}
\newproclaim{assumptionB}{Assumption}
\newproclaim{assumptionK}{Assumption}

\newcommand{\var}{\operatorname{var}}
\newcommand{\E}{\mathrm{E}}
\newcommand{\cov}{\operatorname{cov}}
\newcommand{\ls}{\leq}
\newcommand{\gs}{\geq}
\newcommand{\eps}{\varepsilon}

\newcommand{\dto}{\stackrel{\mathcal{L}}{\longrightarrow}}
\newcommand{\ddeq}{\stackrel{\mathcal{L}^*}{=}}
\newcommand{\pto}{\stackrel{P}{\longrightarrow}}
\newcommand{\deq}{\stackrel{\mathcal{L}}{=}}

\newcommand{\wv}{\widehat{V}}

\makeatother

\begin{document}
\begin{frontmatter}

\title{TFT-bootstrap: Resampling time series in the frequency domain
to obtain replicates in~the~time~domain\thanksref{T1}}
\runtitle{TFT-bootstrap}

\thankstext{T1}{Supported by the DFG graduate college ``Mathematik und Praxis.''}

\begin{aug}
\author[A]{\fnms{Claudia} \snm{Kirch}\thanksref{t1}\ead[label=e1]{claudia.kirch@kit.edu}}
and
\author[B]{\fnms{Dimitris N.} \snm{Politis}\corref{}\thanksref{t2}\ead[label=e2]{dpolitis@ucsd.edu}}
\runauthor{C. Kirch and D. N. Politis}
\affiliation{Karlsruhe Institute of Technology \textup{(}KIT\textup{)} and
University~of~California,~San~Diego}
\address[A]{Institute for Stochastics\\
Karlsruhe Institute of Technology (KIT)\\
Kaiserstr. 89\\
D-76133 Karlsruhe\\
Germany\\
\printead{e1}}
\address[B]{Department of Mathematics\\
University of California, San Diego\\
La Jolla, California 92093-0112\\
USA\\
\printead{e2}}
\end{aug}

\thankstext{t1}{Supported in part by the Stifterverband f\"{u}r die
Deutsche Wissenschaft by funds of the Claussen--Simon-trust.}
\thankstext{t2}{Supported in part by NSF Grant DMS-07-06732.}

\received{\smonth{5} \syear{2010}}
\revised{\smonth{10} \syear{2010}}

%
\begin{abstract}
A new time series bootstrap scheme, the time frequency toggle
(TFT)-bootstrap, is proposed. Its basic idea is to bootstrap the
Fourier coefficients of the observed time series, and then to
back-transform them to obtain a bootstrap sample in the time domain.
Related previous proposals, such as the ``surrogate data'' approach,
resampled only the phase of the Fourier coefficients and thus had only
limited validity. By contrast, we show that the appropriate resampling
of phase \textit{and} magnitude, in addition to some smoothing of
Fourier coefficients, yields a bootstrap scheme that mimics the correct
second-order moment structure for a large class of time series
processes. As a main result we obtain a functional limit theorem for
the TFT-bootstrap under a variety of popular ways of frequency domain
bootstrapping. Possible applications of the TFT-bootstrap naturally
arise in change-point analysis and unit-root testing where statistics
are frequently based on functionals of partial sums. Finally, a small
simulation study explores the potential of the TFT-bootstrap for small
samples showing that for the discussed tests in change-point analysis
as well as unit-root testing, it yields better results than the
corresponding asymptotic tests if measured by size and power.
\end{abstract}

%

\end{frontmatter}

\section{Introduction}

Following Efron's seminal paper~\cite{efron79} on the i.i.d. bootstrap,
researchers have been able to apply resampling ideas in
a variety of non-i.i.d.\vadjust{\goodbreak} situations including the interesting
case of dependent data. B\"{u}hlmann~\cite{buehlmann02},
Lahiri~\cite{lahiri03} and Politis~\cite{politis03a} give reviews of the
state-of-the-art in resampling time series and dependent data.

In the last two decades, in particular,
resampling methods in the frequency domain have become increasingly
popular (see Paparoditis~\cite{papa02}
for a recent survey). One of the first papers to that effect was
Franke and H\"{a}rdle~\cite{franke} who proposed a bootstrap
method based on resampling the periodogram in order
to devise confidence intervals for the spectral density. The
idea behind that approach is that a random vector of the periodogram
ordinates
at finitely many frequencies is approximately independent and
exponentially distributed (cf., e.g., Brockwell and
Davis~\cite{brockwelldavis}, Theorem 10.3.1). Later this approach was also
pursued for different set-ups, for example, for ratio statistics such as
autocorrelations by Dahlhaus and Janas~\cite{dahlhausjanas} or in
regression models by Hidalgo~\cite{hidalgo}.
Dahlhaus and Janas~\cite{dahlhausjanas2} suggested a modification of
the periodogram bootstrap which leads to a correct approximation for a
wider class of statistics such as the sample auto-covariance which---in
contrast to the sample autocorrelation---is not a ratio statistic.
Kreiss and Paparoditis~\cite{papakreiss03} propose the
autoregressive-aided periodogram bootstrap where a parametric time
domain bootstrap is combined with a nonparametric frequency domain bootstrap
in order to widen the class of statistics for which the bootstrap is valid.

We will refer to the above methods as periodogram bootstrapping as
all of the statistics of interest there were functionals of the periodogram.
Since these bootstrap methods resample the periodogram, they
generally do not produce bootstrap pseudo-series in the time
domain. A recent exception is a ``hybrid'' bootstrap of Jentsch and
Kreiss~\cite{jentschkreiss10},
that is, an extension of the aforementioned method of Kreiss and
Paparoditis~\cite{papakreiss03}. 

We now wish to focus on two well-known proposals on frequency-domain
bootstrap methods that
also yield replicates in the time domain, notably:

\begin{itemize}
\item The early preprint by Hurvich and Zeger~\cite{hurvichzeger}
who proposed a parametric bootstrap very similar to our TFT wild
bootstrap of Section~\ref{section_descr_boot}, as well as
a nonparametric frequency-domain bootstrap
based on prewhitening via an estimate of
the MA($\infty$) transfer function.
Although never published, this paper has had substantial influence on
time series literature as it helped inspire many of the above
periodogram bootstrap methods.
Note Hurvich and Zeger~\cite{hurvichzeger}
provide some simulations but give
no theoretical justification for their proposed procedures;
indeed, the first theoretical justification for these
ideas is given in the paper at hand as special cases of the
TFT-bootstrap.

\item The ``\textit{surrogate data}'' approach of Theiler et al. \cite
{theileretal92} has received significant attention in the physics
literature. The idea of
the surrogate data method is to bootstrap the phase of the
Fourier coefficients but keep their magnitude
unchanged.
While most of the literature focuses on heuristics and applications,
some mathematical proofs have been recently provided (see Braun and
Kulperger~\cite{braun}, Chan~\cite{chan97},
Mammen and Nandi~\cite{mammennandipp}, and the recent survey by
Maiwald et al.~\cite{mammen08}).
The surrogate data method was developed for the specific purpose of
testing the null hypothesis of time series linearity
and is not applicable in more general settings.
To see why, note that every surrogate sample has exactly the same
periodogram (and mean) as the original sequence.
Hence, the method fails to approximate the distribution of
any statistic that is a function of first- and second-order moments,
thus excluding all cases where periodogram resampling
has proven to be useful; see our Proposition~\ref{prop1} in
Section~\ref{section_diff_methods}.
\end{itemize}

In the paper at hand, we propose to resample the Fourier
coefficients---which can effectively be computed using a fast Fourier
transform (FFT)---in a variety of ways similar to modern periodogram bootstrap
methods, and then obtain time series resamples using an inverse
FFT. Since we start out with an observation sequence in the time
domain, then jump to the frequency domain for resampling just to get
back to the time domain again, we call this type of resampling a
\textit{time frequency toggle} (TFT) bootstrap. The TFT-bootstrap is an
extension of existing periodogram bootstrap methods as it yields almost
identical procedures when applied to statistics based on periodograms,
but it is also applicable in situations where the statistics of
interest are not expressible by periodograms; for more details we refer
to Section~\ref{section_appl}.

The TFT-bootstrap is related to the surrogate data approach but is more
general since it also resamples the magnitudes of Fourier coefficients
and not just their phases. As a result, the TFT is able to correctly
capture the distribution of statistics that are based on the periodogram.
The TFT, however, shares with the surrogate data approach the
inability to approximate the distribution of the sample mean; luckily,
there are plenty of methods in the bootstrap literature to accomplish
that, for example, the block bootstrap and its variations, the
AR-sieve bootstrap, etc. (for details, see Lahiri~\cite{lahiri03}, B\"
{u}hlmann~\cite{buehlmann02},
Politis~\cite{politis03a}).

In this paper we provide some general theory for the TFT-bootstrap
which not only gives a long-due theoretical justification for one of
the proposals by Hurvich and Zeger~\cite{hurvichzeger} but also allows
for several modern extensions of these early ideas.
In particular, we prove that the TFT sample has asymptotically the
correct second-order moment structure (Lemma~\ref{lem_cov}) and
provide a functional central limit theorem (FCLT, Theorem \ref
{th_main} and Corollary~\ref{cor_main}) for the TFT-sample. This is a
much stronger result than the asymptotic normality with correct
covariance structure of a finite subset as proved, for example, by
Braun and Kulperger~\cite{braun} for the
surrogate data method.

As in the surrogate data method, the TFT sample paths are shown to be
(asymptotically) Gaussian; so in a sense the TFT
approximates possibly nonlinear time\vadjust{\goodbreak} series with a Gaussian process
having the correct second-order moment structure. This seems to be
inevitable in all methods using discrete Fourier transforms due to the
fact that Fourier coefficients are asymptotically normal under very
general assumptions.
However, in contrast to the surrogate data method, the TFT is able to
capture the distribution of many useful statistics (cf. Section \ref
{section_appl}).
For example, our FCLT implies the validity of inference for statistics
such as CUSUM-type statistics in change-point analysis (cf.
Section~\ref{section_cpa}) or least-squares statistics in unit-root
testing (cf. Section~\ref{section_unitroot}). The TFT-bootstrap is
also valid
for periodogram-based (ratio) statistics such as sample
autocorrelations or Yule--Walker estimators; this validity is
inherited by the corresponding results of the periodogram bootstrapping
employed for the TFT (cf. Section~\ref{sec_stat_perio}).

Furthermore, in many practical situations one does not directly observe
a stationary sequence but needs to estimate it first. In Corollary \ref
{cor_boot_est} we prove the validity of the TFT-bootstrap when applied
to such estimated sequences. For example, in change-point analysis
(Section~\ref{section_cpa}) as well as unit-root testing (Section \ref
{section_unitroot}) one can use estimators to obtain an approximation
of the underlying stationary sequence under the null hypothesis as well
as under the alternative. As in both examples the null hypothesis is
that of a stationary sequence; this feature enables us to construct
bootstrap tests that capture the null distribution of the statistic in question
even when presented with data that obey the alternative hypothesis. As
a consequence, these bootstrap tests asymptotically capture the correct
critical value even under the alternative hypothesis which not only
leads to the correct size of the tests but also to a good power behavior.

The remainder of the paper is organized as follows.
In the next section we give a detailed description on how the
TFT-bootstrap works. In particular we describe several specific
possibilities of how to get pseudo-Fourier coefficients. In
Section~\ref{section_fclt} we state the main theorem, a functional
limit theorem for the TFT-bootstrap. The FCLT holds true under certain
high-level assumptions on the bootstrapped Fourier coefficients; these
are further explored in Sections~\ref{section_valboot} and \ref
{section_prop_freq}. In particular it is shown that the TFT-bootstrap
replicates the correct second-order moment structure for a large class
of observed processes including nonlinear processes (cf. Section \ref
{section_prop_freq}).
Finally, we prove the validity of the TFT-bootstrap for certain
applications such as unit-root testing or change-point tests in
Section~\ref{section_appl} and explore the small sample performance in
the simulation study of Section~\ref{section_sim}.
Our conclusions are summarized in Section~\ref{section_conclusions}.
Proofs are sketched in Section~\ref{section_fclt_proof}, while the
complete technical proofs can be found in electronic supplementary material.

\section{Description of the TFT-bootstrap}\label{section_descr_boot}

Assume we have observed $V(1),\ldots,\break V(T)$, where
\renewcommand{\theassumptionP}{$\mathcal{P}$.\arabic{assumptionP}}
\begin{assumptionP}
\label{ass_W}
$\{V(i)\dvtx i\gs1\}$ is a stationary process with absolutely summable
auto-covariance function $\gamma(\cdot)$. In this\vadjust{\goodbreak} case the spectral
density of the process exists, is continuous and bounded. It is defined by
%
\begin{equation}\label{eq_pres_spec_dens}
f(\lambda)=\frac{1}{2\pi}\sum_{n=-\infty}^{\infty}e^{-in\lambda
}\gamma(n)
\end{equation}
(see, e.g., Brockwell and Davis~\cite{brockwelldavis}, Corollary 4.3.2).
\end{assumptionP}

Since we will prove a functional central limit theorem for the
bootstrap sequence, the procedure only makes sense if the original
process fulfills the same limit theorem.
\begin{assumptionP}\label{ass_clt}
$\{V(i)\dvtx i\gs1\}$ fulfills the following functional central limit theorem:
\[
\Biggl\{\frac{1}{\sqrt{2\pi f(0) T}}\sum_{t=1}^{\lfloor T u\rfloor
}\bigl(V(t)-\E V(1)\bigr)\dvtx0\ls u\ls1\Biggr\}\dto\{W(u)\dvtx0\ls u\ls1\},
\]
where $f(\cdot)$ is the spectral density of $\{V(\cdot)\}$ and $\{
W(\cdot)\}$ is a standard Wiener process.
\end{assumptionP}

We may need the following assumption on the spectral density:
\begin{assumptionP}\label{ass_process_dens}
Let the spectral density be bounded from below $f(\lambda)\gs c>0$ for
all $0\ls\lambda\ls2\pi$.
\end{assumptionP}

We denote by $Z(j)=V(j)-\bar{V}_T$ the centered counterpart of the observations
where $\bar{V}_T= T^{-1}\sum_t V(t)$.
Consider the FFT coefficients of the observed stretch, that is,
\begin{eqnarray*}
x(j)&=&\frac{1}{\sqrt{T}}\sum_{t=1}^TV(t)\cos(-\lambda_j t),\\
y(j)&=&\frac{1}{\sqrt{T}}\sum_{t=1}^TV(t)\sin(-\lambda_j t),
\end{eqnarray*}
thus
\[
x(j)+iy(j)=\frac{1}{\sqrt{T}}\sum
_{t=1}^TV(t)\exp(-
i \lambda_j t),
\]
where $\lambda_j=2\pi j/T$ for $j=1,\ldots,T$. Note that the Fourier
coefficients $x(j)$, $y(j)$ depend on $T$, but to keep the notation
simple we suppress this dependence.

The principal idea behind all bootstrap methods in the frequency domain
is to make use of the fact that the Fourier coefficients
$x(1),y(1),\ldots,x(N),\break y(N)$
are asymptotically independent\vadjust{\goodbreak} and normally distributed, where
$N=\lfloor(T-1)/2\rfloor$ denotes the largest integer smaller or
equal to $(T-1)/2$,
and
%
\begin{eqnarray}\label{eq_fft_distr}
\E x(j)&\to&0,\qquad \E y(j)\to0 \quad\mbox{and}\nonumber\\[-8pt]\\[-8pt]
\var x(j)&=&\pi f(\lambda_j)+o(1),\qquad \var y(j)= \pi
f(\lambda_j)+o(1)\nonumber
\end{eqnarray}
for $j=1,\ldots,N$
as $T\to\infty$ where $f(\cdot)$ is the spectral density (see, e.g.,
Chapter 4
of Brillinger~\cite{brillinger} for a precise formulation of this
vague statement). Lahiri~\cite{lahiri03a} gives necessary as well as
sufficient conditions for the asymptotic independence and normality of
tapered as well as nontapered Fourier coefficients for a much larger
class of time series not limited to linear processes in the strict
sense. Shao and Wu~\cite{shaowu07} prove this statement uniformly over
all finite subsets for a large class of linear as well as nonlinear
processes with nonvanishing spectral density. The uniformity of their
result is very helpful since it implies convergence of the
corresponding empirical distribution function (see the proof of
Lemma~\ref{lem_edf} below).

The better known result on the periodogram
ordinates states that $x(j)^2+y(j)^2$ are asymptotic independent
exponentially distributed with expectation $2\pi f(\lambda_j)$ (see,
e.g., Chapter 10 of Brockwell and Davis~\cite{brockwelldavis}).
The latter is what most bootstrap versions are based on.
By contrast, our TFT-bootstrap will focus on the former property
(\ref{eq_fft_distr}), that is, the fact that $x(j),y(j)$ are
asymptotically i.i.d. $N(0,\pi f(\lambda_j))$.

Let us now recall some structural properties of the Fourier
coefficients, which are important in order to understand the procedure
below. First note that
%
\begin{equation}\label{eq_conj_com}
x(T-j)=x(j)  \quad\mbox{and}\quad  y(T-j)=-y(j).
\end{equation}
This symmetry implies that
the Fourier coefficients for $j=1,\ldots, N$ carry all the necessary
information required in order to
recapture (by inverse FFT) the original series up to
an additive constant.

The symmetry relation (\ref{eq_conj_com}) shows in particular that all the
information carried in the coefficients for
$\lambda_j$, $j>T/2$, is already contained in the coefficients for
$\lambda_j$, $j\ls N$.
The information that is missing is the information about the mean of
the time series which
is carried by the remaining coefficients belonging to $\lambda_T$ and
$\lambda_{T/2}$ (the latter only
when $T$ is even).

To elaborate,
\[
x(T)=\sqrt{T}\bar{V}_T
\]
carries the information about the mean of the observations $V(1),\ldots
,V(T)$; moreover,
\[
y(T)=0.
\]
For $T$ even, we further have some additional information about the
``alternating'' mean
\[
y(T/2)=0 \quad\mbox{and}\quad  x(T/2)=\frac{1}{\sqrt{T}} \sum
_{t=1}^T(-1)^tV(t).\vadjust{\goodbreak}
\]
Note that the value of
the FFT coefficients for $j=1,\ldots, N$ is the same for a sequence
$V(t)-c$ for all $c\in\mathbb{R}$.
Hence those Fourier coefficients are invariant under additive constants
and thus contain no information about the mean.
Similarly all the information about the bootstrap mean is carried only
in the bootstrap version of $x(T)$ [as well as $x(T/2)$]. The problem
of bootstrapping the mean is therefore separated from getting a time
series with the appropriate covariance structure and will not be
considered here. In fact, we show that any asymptotically correct
bootstrap of the mean if added to our bootstrap time series [cf. (\ref
{eq_form_bs2})] yields the same asymptotic behavior in terms of its
partial sums as the original uncentered time series.

Our procedure works as follows:
\begin{longlist}[\textit{Step} 1:]
\item[\textit{Step} 1:] Calculate the Fourier coefficients using the fast Fourier
transform (FFT) algorithm.
\item[\textit{Step} 2:] Let $y^*(T)=x^*(T)=0$; if $T$ is even, additionally let
$x^*(T/2)=y^*(T/2)=0$.
\item[\textit{Step} 3:] Obtain a bootstrap sequence $x^*(1),y^*(1),\ldots,x^*(N),y^*(N)$
using, for example, one of bootstrap procedures described below.
\item[\textit{Step} 4:] Set the remaining bootstrap Fourier coefficients according to
(\ref{eq_conj_com}), that is, $x^*(T-j)=x^*(j)$ and $y^*(T-j)=-y^*(j)$.
\item[\textit{Step} 5:] Use the inverse FFT algorithm to transform the bootstrap Fourier
coefficients $x^*(j)+iy^*(j)$, $j=1,\ldots,T$, back into the time domain.
\end{longlist}
We
thus obtain a bootstrap sequence $\{Z^*(t)\dvtx1\ls t\ls T\}$ which is
real-valued and centered, and can be used for inference on a large
class of statistics that are based on partial sums of the centered
process $\{Z(\cdot)\}$;
see Section~\ref{section_cpa} for examples.
\begin{rem}
Note that the exact form of $Z^*(t)$ is the following:
%
\begin{eqnarray}\label{eq_form_bs}
Z^*(t)&=&\frac{1}{\sqrt{T}}\sum_{j=1}^T\bigl(x^*(j)+iy^*(j)\bigr)\exp(2\pi i t
j/T)\nonumber\\[-8pt]\\[-8pt]
&=&\frac{2}{\sqrt{T}}\sum_{j=1}^N\bigl( x^*(j)\cos(2\pi t
j/T)-y^*(j)\sin(2\pi tj/T) \bigr).\nonumber
\end{eqnarray}
\end{rem}
\begin{rem}\label{rem_centering}
In order to obtain a bootstrap sequence of the noncentered observation
process $\{V(\cdot)\}$ we can add a bootstrap mean $\mu_T^{*}$ to the
$\{Z^*(\cdot)\}$
process; here, $\mu_T^{*}$
is obtained by a separate bootstrap process independently from $ \{
Z^{*}(\cdot)\}$, which is asymptotically normal with the correct
variance, that is, it fulfills (\ref{eq_boot_mean}).

Precisely, the bootstrap sequence
%
\begin{equation}
\label{eq_form_bs2}
V^{*}(t)=Z^{*}(t)+\mu_T^{*}
\end{equation}
gives a bootstrap approximation of $\{V(\cdot)\}$. Here, $\{Z^*(\cdot
)\}$ contains the information about the covariance structure of the
time series, and $\mu^*$ contains the information of the sample mean
as a random variable of the time series. How to obtain the latter will
not be considered in this paper.
\end{rem}

In Corollary~\ref{cor_boot_est} we give some conditions under which
the above procedure remains asymptotically valid if instead of the
process $V(\cdot)$ we use an estimated process $\widehat{V}(\cdot)$;
this is important in some applications.

Now we are ready to state some popular bootstrap algorithms in the
frequency domain. We have adapted them in order to bootstrap the
Fourier coefficients rather than the periodograms. Our procedure can
easily be extended to different approaches.

\subsection*{Residual-based bootstrap (RB)}

\mbox{}

\textit{Step} 1: First estimate the spectral density $f(\cdot
)$ by $\widehat{f}(\cdot)$ satisfying
%
\begin{equation}\label{eq_est_spec_dens}
{\sup_{\lambda\in[0,\pi]}}|\widehat{f}(\lambda)-f(\lambda
)|\pto0.
\end{equation}
This will be denoted Assumption~\ref{ass_a1} in Section
\ref{section_prop_freq}.
Robinson~\cite{robinson91} proves such a result for certain kernel
estimates of the spectral density based on periodograms for a large
class of processes including but not limited to linear processes. For
linear processes he also proves the consistency of the spectral density
estimate as given above when an automatic bandwidth selection procedure
is used.
Shao and Wu~\cite{shaowu07} also prove this result for certain kernel
estimates of the spectral density for processes satisfying some
geometric-moment contraction condition, which includes a large class of
nonlinear processes. Both results are summarized in Lem\-ma~\ref
{lem_spectralestimate}.

\textit{Step} 2: Next estimate the residuals of the real, as
well as imaginary, part of the Fourier coefficients and put them
together into a vector $\{\widetilde{s}_j\dvtx 1\ls j\ls2N\}$; precisely
let
\[
\widetilde{s}_j=
\frac{x(j)}{\sqrt{\pi\widehat{f}(\lambda_j)}}, \qquad \widetilde{s}_{N+j}=
\frac{y(j)}{\sqrt{\pi\widehat{f}(\lambda_{j})}},
\]
$j=1,\ldots, N$. Then standardize them, that is, let
\[
s_j=\frac{\widetilde{s}_j-{1}/({2N})\sum_{l=1}^{2N}\widetilde
{s}_l}{\sqrt{{1}/({2N})\sum_{t=1}^{2N}(\widetilde
{s}_t-{1}/({2N})\sum_{l=1}^{2N}\widetilde{s}_l)^2}}.
\]
Heuristically these residuals are approximately i.i.d., so that i.i.d.
resampling methods are reasonable.

\textit{Step} 3: Let $s_j^*$, $j=1,\ldots,2N$, denote an i.i.d.
sample drawn randomly
and with replacement from $s_1,\ldots,s_{2N}$. As usual, the
resampling step is
performed conditionally on the data $V(1),\ldots, V(T)$.

\textit{Step} 4: Define the bootstrapped Fourier coefficients by
%
\begin{equation}\label{eq_bc_b1}
x^*(j)=\sqrt{\pi\widehat{f}(\lambda_j)}s_j^*, \qquad
y^*(j)=\sqrt{\pi\widehat{f}(\lambda_j)}s_{N+j}^* .
\end{equation}
An analogous approach---albeit focusing on the periodogram
ordinates instead of the FFT coefficients---was proposed by Franke and
H\"{a}rdle~\cite{franke} in order to yield a bootstrap distribution of
kernel spectral density
estimators.

\subsection*{Wild bootstrap (WB)}

The wild bootstrap also makes use of an estimated spectral
density further exploiting the knowledge about the asymptotic normal
distribution of the Fourier coefficients. Precisely, the
WB replaces $s_j^*$ above by independent standard normal distributed
random variables $\{G_j\dvtx 1\ls j\ls2N\}$ in order to obtain the
bootstrap Fourier coefficients as in (\ref{eq_bc_b1}).
This bootstrap was already suggested by Hurvich and Zeger \cite
{hurvichzeger}, who considered it in a simulation study, but did not
obtain any theoretical results.

An analogous approach---albeit focusing on the
periodogram---was discussed by Franke and H\"{a}rdle~\cite{franke} who
proposed multiplying the periodogram with i.i.d. exponential random variables.

\subsection*{Local bootstrap (LB)}

The advantage of the local bootstrap is that it does not need
an initial estimation of the spectral density. The idea is that in a
neighborhood of each frequency the distribution of the different
coefficients is almost identical (if the spectral density is smooth).
It might therefore be better able to preserve some information beyond
the spectral density that is contained in the Fourier coefficients.
An analogous procedure for periodogram ordinates was first proposed by
Paparoditis and Politis~\cite{papapolitis99}.
For the sake of simplicity we will only consider bootstrap schemes that
are related to kernels.

Recall that $x(-j)=x(j)$, $ x(\lceil T/2\rceil+j)=x(N-j)$ and
$y(-j)=-y(j)$, $y(\lceil T/2\rceil+j)=y(N-j)$ for $j=1,\ldots,N+1$,
and, for $T$ even, $x(N+1)=y(N+1)=\frac{1}{\sqrt{T}} \sum
_{t=1}^T(-1)^tV(t)$ ($\lceil x \rceil$ denotes the smallest integer
larger or equal than $x$).
Furthermore let $x(0)=y(0)=0$ be the Fourier coefficients of the
centered sequence. For $j<-T/2$ and $j>T/2$ the coefficients are
periodically extended with period $T$.\vspace*{8pt}

\textit{Step} 1: Select a symmetric, nonnegative kernel
$K(\cdot)$ with
\[
\int K(t) \,dt=1.
\]
In Section~\ref{section_prop_freq} we assume some additional
regularity conditions on the kernel in order to get the desired results.
Moreover select a bandwidth $h_T$ fulfilling $h_T\to0$ but $Th_T\to
\infty$.

\textit{Step} 2: Define i.i.d. random variables $J_{1,T},\ldots
, J_{2N,T}$ on $\mathbb{Z}$ with
%
\begin{equation}\label{eq_weights_local}
p_{s,T}=P(J_{j,T}=s)=\frac{K( {2\pi s}/({T h_T})
)}{\sum_{j=-\infty}^{\infty}K( {2\pi j}/({T h_T}) )}.
\end{equation}
Independent of them define i.i.d. Bernoulli r.v. $B_1,\ldots,B_{2N}$
with parameter $1/2$.

\textit{Step} 3: Consider now the following bootstrap sample:
\[
\widetilde{x}^*(j)=\cases{
x(j+J_{j,T}), &\quad if $B_j=0$,\cr
y(j+J_{j,T}), &\quad if $B_j=1$,}
\]
and
\[
\widetilde{y}^*(j)=\cases{
y(j+J_{N+j,T}), &\quad if $B_{N+j}=0$,\cr
x(j+J_{N+j,T}), &\quad if $B_{N+j}=1$.}
\]
Finally the bootstrap Fourier coefficients are defined as the centered
versions of~$\widetilde{x}^*$, respectively, $\widetilde{y}^*$, namely by
\begin{eqnarray*}
x^*(j)&=&\widetilde{x}^*(j)-\frac{1}{2}\sum_{s\in\mathbb
{Z}}p_{s,T}\bigl( x(j+s)+y(j+s) \bigr),\\
y^*(j)&=&\widetilde{y}^*(j)-\frac{1}{2}\sum_{s\in\mathbb
{Z}}p_{s,T}\bigl( x(j+s)+y(j+s) \bigr).
\end{eqnarray*}
This is slightly different from Paparoditis and Politis \cite
{papapolitis99}, since they require that $x^*(j)$ and $y^*(j)$ share
the same $J_{j,T}$ which is reasonable if one is interested in
bootstrapping the periodogram but not necessary for bootstrapping the
Fourier coefficients.


\subsection*{Comparison of the three methods}

The aforementioned three bootstrap methods,
residual bootstrap (RB), wild bootstrap (WB) and local bootstrap (LB),
are all first-order consistent under standard conditions. A rigorous
theoretical comparison would entail higher-order considerations which
are not available in the literature and are beyond the scope of this work.
Intuitively, one would expect the RB and LB procedures to perform
similarly in applications since these two bootstap methods
share a common underlying idea, that is, that nearby periodogram/FFT
ordinates are i.i.d. By contrast, the WB involves the
generation of extraneous Gaussian random variables
thus forcing the time-domain bootstrap sample paths to be
Gaussian. For this reason alone, it is
expected that if a higher-order property holds true
in our setting, it will likely be shared by RB and LB but not WB. Our
finite-sample simulations in Section~\ref{section_sim} may hopefully
shed some
additional light on the comparison between RB and LB.

\subsection{Comparison with other frequency domain methods}
\label{section_diff_methods}
First, note that the TFT wild bootstrap is identical to the
parametric frequency-domain
bootstrap proposal of Hurvich and\vadjust{\goodbreak} Zeger~\cite{hurvichzeger}.
By contrast, the nonparametric bootstrap
proposal of Hurvich and Zeger~\cite{hurvichzeger} was
based on prewhitening via an estimate of
the MA($\infty$) transfer function. Estimating the
transfer function presents an undesirable complication since
prewhitening can be done in an easier fashion
using any consistent estimator of the
spectral density; the residual-based TFT
exploits this idea based on the work of
Franke and H\"{a}rdle~\cite{franke}.
The local bootstrap TFT is a more modern extension of the
same underlying principle, that is, exploiting the approximate
independence (but not i.i.d.-ness) of periodogram ordinates.

We now attempt to shed some light on the relation between the TFT and
the surrogate data method of Theiler et al.~\cite{theileretal92}.
Recall that the surrogate data approach amounts to using
%
\begin{equation}\label{eq.T1234}
\sqrt{I(j)}\cos(2\pi U_j)+i\sqrt{I(j)}\sin(2\pi U_j)
\end{equation}
as bootstrap Fourier coefficients at point $j$
where $I(j)={x^2(j)+y^2(j)}$ is the periodogram at point $j$, and $\{
U_j\dvtx j\}$ are i.i.d. uniform on $[0,1]$.
So the periodogram (and mean) computed from surrogate is identical
to the original ones.
As a result we have the following proposition.
\begin{prop}\label{prop1}
The surrogate data method fails to approximate the distribution
of any statistic that can be written as a function of the periodogram
(and/or sample mean) as long as this statistic has nonzero
large-sample variance.
\end{prop}

The proof of this proposition is obvious since by replicating
the periodogram (and sample mean) exactly, the surrogate approach will
by necessity approximate the distribution of the statistics in question
by a point mass, that is, zero variance.
Hence, the surrogate data method will not be able to correctly capture
distributional properties for a large class of statistics
that are based entirely on periodograms and the mean (cf. also
Chan~\cite{chan97}). In hypothesis testing, this failure would result in
having power equal to the size when the surrogate approach is applied
to obtain critical values for test statistics based entirely on
periodogram and/or mean; at the very least, some loss of power in other
situations can be expected.

The aforementioned TFT methods obviously do not have this disadvantage;
in fact, they can be successfully applied to this precise class of statistics
(cf. Section~\ref{sec_stat_perio}). 

For comparison purposes, we now describe a
\textit{nonsmoothed} wild bootstrap that is the closest relative of the
surrogate data method that fits in our framework.
Note that all TFT-bootstrap schemes involve smoothing in the frequency
domain before resampling (cf. also Assumption \ref
{ass_boot_2b}); however, one could consider bootstrapping without this
smoothing step.
As the wild bootstrap works by multiplying normal random variables with
an estimator of the spectral density, that is, a smoothed version of
the periodogram, the nonsmoothed wild bootstrap multiplies
normal random variables with the original Fourier coefficients
$x(j),y(j)$. By the Box--Muller transform, the nonsmoothed wild
bootstrap gives
the following as the bootstrap (complex-valued) Fourier coefficient at
point $j$:
%
\begin{equation}\label{eq.S1234}
x(j)\sqrt{-2\log(\widetilde{U}_j)}\cos(2\pi U_j)+iy(j)\sqrt{-2\log
(\widetilde{U}_j)}\sin(2\pi U_j),
\end{equation}
where $\{U_j\dvtx j\},\{\widetilde{U}_j\dvtx j\}$ are i.i.d. uniform on $[0,1]$
independent from each other.

Comparing equation (\ref{eq.S1234}) to equation (\ref{eq.T1234}) we
see that the surrogate data approach is closely related to the
nonsmoothed wild bootstrap; the main difference is that the wild
bootstrap does not only bootstrap the phase but also the magnitude of
the Fourier coefficients. Nevertheless, the nonsmoothed wild bootstrap
does not suffer from the severe deficiency outlined in Proposition~\ref{prop1}
since it does manage to
capture the variability of the periodogram to some extent.
To elaborate, note that
it is possible to prove a functional limit theorem
[like our Theorem~\ref{th_main}(a) in the next section]
for the nonsmoothed wild bootstrap but only under the provision that
a \textit{smaller resample size} is employed;
that is, only a fraction of the bootstrap sample is used
to construct the partial sum process ($m/T\to0$). This
undersampling condition is necessary here since without it the
asymptotic covariance structure would not be correct.
Hence, even the nonsmoothed wild bootstrap, although crude, seems
(a) preferable to the surrogate data method and (b)
inferior with respect to the TFT-bootstrap; this relative
performance comparison is clearly born out in simulations
that are not reported here due to lack of space.

\section{Functional limit theorem for the bootstrap sample}\label
{section_fclt}

In this section we state the main result, namely a functional limit
theorem for the partial sum processes of the bootstrap sample.

The theorem is formulated in a general way under some meta-assumptions
on the resampling scheme in the frequency domain that ensure the
functional limit theorem back in the time domain. In Section \ref
{section_valboot} we verify those conditions for the bootstrap schemes
given in the previous section. We would like to point out that the
meta-assumptions we give are the analogues of what is usually proved
for the corresponding resampling schemes of the periodograms, which are
known to hold for a large class of processes.

The usage of meta-assumptions allows the reader to extend results to
different bootstrap schemes in the frequency domain.

By $\E^*$, $\var^*$, $\cov^*$ and $P^*$ we denote as usual the
bootstrap expectation, variance, covariance and probability. 
We essentially investigate three sets of assumptions.

The first one is already implied by the above mentioned bootstrap schemes.

\renewcommand{\theassumptionB}{$\mathcal{B}$.\arabic{assumptionB}}
\begin{assumptionB}\label{ass_boot_1}
For the bootstrap scheme in the frequency domain, the coefficients
$\{x^*(k)\dvtx 1\ls k\ls N\}$ and $\{y^*(k)\dvtx 1\ls k\ls N\}$ are independent
sequences as well as mutually independent (conditionally on the data) with
\[
\E^*(x^*(k))=\E^*(y^*(k))=0.
\]
\end{assumptionB}
\begin{rem}\label{rem_moments}
Instead of assuming that the bootstrap samples are already
centered it is sufficient that the bootstrap means in the frequency
domain converge uniformly to $0$ with a certain rate, that is,
\[
\sup_{1\ls k\ls N}\bigl(|\E^*(x^*(k))| +|\E
^*(y^*(k))| \bigr)=o_P\Biggl( \sqrt{\frac{m}{T\log^2 T}}
\Biggr),
\]
where $m$ is the parameter figuring in Lemma~\ref{lem_cov} (resp.,
Theorem~\ref{th_main} below).
\end{rem}

%
\begin{assumptionB}\label{ass_boot_2b}
Uniform convergence of the second moments of the bootstrap sequence in
the frequency domain, that is,
\begin{eqnarray*}
{\sup_{1\ls k\ls N}}|\var^*(x^*(k))-\pi f(\lambda_k)
|&=&o_P(1),\\
{\sup_{1\ls k\ls N}}|\var^*(y^*(k))-\pi f(\lambda
_k)|&=&o_P(1).
\end{eqnarray*}
\end{assumptionB}
\begin{assumptionB}\label{ass_boot_4}
Uniform boundedness of the fourth moments of the bootstrap sequence in
the frequency domain
\[
\sup_{1\ls k\ls N}\E^*(x^*(k))^4\ls C +o_P(1), \qquad\sup_{1\ls k\ls
N}\E^*(y^*(k))^4\ls C +o_P(1).
\]
\end{assumptionB}

Let us now recall the definition of the Mallows distance on the space
of all real Borel probability measures with finite variance. It is
defined as
\[
d_2(P_1,P_2)=\inf( \E|X_1-X_2|^2 )^{1/2},
\]
where the infimum is taken over all real-valued variables $(X_1,X_2)$
with marginal distributions $P_1$ and $P_2$, respectively. Mallows \cite
{mallows72} has proved the equivalence of convergence in this metric
with distributional convergence in addition to convergence of the
second moments. The results remain true if we have convergence in a
uniform way as in Assumption~\ref{ass_boot_3} below.
This shows that Assumption~\ref{ass_boot_3} implies
Assumption~\ref{ass_boot_2b}.
\begin{assumptionB}\label{ass_boot_3}
Let the bootstrap scheme in the frequency domain converge uniformly in
the Mallows distance to the same limit as the Fourier coefficients do
\begin{eqnarray*}
\sup_{1\ls j\ls N}d_2(\mathcal{L}^*(x^*(j)),N(0,\pi f(\lambda
_j)))&=&o_P(1), \\
\sup_{1\ls j\ls N}d_2(\mathcal{L}^*(y^*(j)),N(0,\pi f(\lambda_j)))&=&o_P(1).
\end{eqnarray*}
\end{assumptionB}

We will start with some results concerning the asymptotic covariance
structure of the partial sum process; all asymptotic results are taken
as $T \to\infty$.
\begin{lemma}\label{lem_cov}
\begin{longlist}[(a)]
\item[(a)] Let Assumption~\ref{ass_boot_1} be fulfilled.
Then, for any $0\ls u\ls1$ and \mbox{$1\ls m\ls T$},
\[
\E^*\Biggl( \frac{1}{\sqrt{m}}\sum_{l=1}^{\lfloor m u\rfloor
}Z^*(l) \Biggr)=0.
\]
\item[(b)] Let Assumptions~\ref{ass_W}, \ref
{ass_boot_1} and~\ref{ass_boot_2b} be fulfilled. Then,
for $0\ls u,v\ls1$,
\begin{eqnarray*}
&&\cov^*\Biggl( \frac{1}{\sqrt{m}}\sum_{l_1=1}^{\lfloor m u\rfloor
}Z^*(l_1), \frac{1}{\sqrt{m}}\sum_{l_2=1}^{\lfloor m v\rfloor
}Z^*(l_2)\Biggr)\\
&&\qquad\pto
\cases{2\pi f(0) \min(u,v), &\quad $\dfrac m T \to0$,\vspace*{2pt}\cr
2\pi f(0)[\min(u,v)-uv],&\quad $m=T$.}
\end{eqnarray*}
\item[(c)] Moreover, under Assumptions~\ref{ass_W}, \ref
{ass_boot_1} and~\ref{ass_boot_2b},
\[
\cov^*(Z^*(l_1),Z^*(l_2))=\cov(V(l_1),V(l_2))+o_P(1)=\cov
(Z(l_1),Z(l_2))+o_P(1)
\]
for all fixed $l_1, l_2$.
\end{longlist}
\end{lemma}

As already pointed, out using frequency domain methods separates the
problem of an appropriate bootstrap mean from the problem of obtaining
a bootstrap sample with the appropriate covariance structure. As a
result, the bootstrap sample $\{Z^*(\cdot)\}$ is centered and thus the
bootstrap version of the centered time series $\{Z(\cdot)\}=\{V(\cdot
)-\bar{V}_T\}$.

The above lemma shows that the bootstrap process $\{Z^*(\cdot)\}$ as
well as its partial sum process has the correct auto-covariance
structure. The following theorem gives a functional central limit
theorem in the bootstrap world, showing that the bootstrap partial sum
process also has the correct second-order moment structure.
In fact, the partial sum process of a centered time series converges to
a Brownian bridge, while the subsampled partial sum processes converges
to a Wiener process. As the following theorem shows this behavior is
exactly mimicked by our TFT-bootstrap sample.
\begin{theorem}\label{th_main}
Let Assumptions~\ref{ass_W}, \ref
{ass_boot_1}--\ref{ass_boot_4} be fulfilled.

\begin{longlist}[(a)]
\item[(a)]
If $m/T\to0$, then it holds (in probability)
\[
\Biggl\{\frac{1}{\sqrt{ 2\pi f(0)m}}\sum_{l=1}^{\lfloor m u\rfloor
}Z^*(l)\dvtx0\ls u\ls1 \Big|  V(\cdot)\Biggr\}\stackrel
{D[0,1]}{\longrightarrow}\{W(u)\dvtx0\ls u\ls1\},
\]
where $\{W(u)\dvtx0\ls u\ls1\}$ is a Wiener process.

\item[(b)]
If additionally Assumption~\ref{ass_boot_3} is
fulfilled, we obtain (in probability)
\[
\Biggl\{\frac{1}{\sqrt{2\pi f(0) T}}\sum_{l=1}^{\lfloor T u\rfloor
}Z^*(l)\dvtx0\ls u\ls1 \Big|  V(\cdot)\Biggr\}\stackrel
{D[0,1]}{\longrightarrow}\{B(u)\dvtx0\ls u\ls1\},
\]
where $\{B(u)\dvtx0\ls u\ls1\}$ is a Brownian bridge.
\end{longlist}
\end{theorem}
\begin{rem}\label{rem_lind}
The stronger Assumption~\ref{ass_boot_3} is needed only
to get asymptotic normality of the partial sum process, that is, part
(b) above. In the proof of Theorem~\ref{th_main} for $m/T\to0$ we use
the Lindeberg condition to obtain asymptotic normality. However, for
$m=T$ the latter is not fulfilled because the variance of single
summands (e.g., $l=T$) is not negligible anymore. But for the same
reason the Feller condition is also not fulfilled which means that we
cannot conclude that the sequence is not asymptotically normal. In
fact, failure of asymptotic normality is hard to imagine in view of
Corollary~\ref{cor_braun}. Therefore we recommend to always use $m=T$
in applications even in situations where Assumption \ref
{ass_boot_3} is hard to verify.
\end{rem}
\begin{rem}\label{rem_surrogate_1}The bootstrap variance is usually
related to the periodogram which is not consistent without smoothing. Therefore
Assumption~\ref{ass_boot_2b} ensures that the bootstrap
scheme includes some smoothing. For the bootstrap, however, this is not
entirely necessary, and we can also bootstrap without smoothing first.
The simplest example is the nonsmoothed wild bootstrap as described in
Section~\ref{section_diff_methods}. One can then still prove the
result of Theorem~\ref{th_main}, but only for $m/T\to0$. In this
situation this condition is necessary, since without it the asymptotic
covariance structure is not correct (i.e., the assertion of Lemma \ref
{lem_cov} is only true for $m/T\to0$), $m\approx\sqrt{T}$ would be a
good rule of thumb. While this is a very simple approach (without any
additional parameters) it does not give as good results as the
procedure we propose. Heuristically, this still works because the
back-transformation does the smoothing, but to obtain a consistent
bootstrap procedure we either need some smoothing in the frequency
domain as in Assumption~\ref{ass_boot_2b} or do some
under-sampling back in the time domain, that is, $m/T\to0$. In fact,
one of the main differences between our wild TFT-bootstrap and the
surrogate data approach is the fact that the latter does not involve
any smoothing in the frequency domain. The other difference being that
the surrogate data approach only resamples the phase but not the
magnitude of the Fourier coefficients. For more details we refer to
Section~\ref{section_diff_methods}.
\end{rem}

Some applications are based on partial sums rather than centered
partial sums. This can be obtained as described in Remark \ref
{rem_centering}. The following corollary then is an immediate
consequence of Theorem~\ref{th_main}(b).
\begin{cor}\label{cor_main}
Let the assumptions of Theorem~\ref{th_main}\textup{(b)} be fulfilled. Let
$\mu_T^*$ be a bootstrap version of the mean $\mu=\E V(1)$ [taken
independently from $\{Z^*(\cdot)\}$] such that for all
$z\in\mathbb{R}$,
%
\begin{equation}\label{eq_boot_mean}
P^{*}\bigl( \sqrt{T}(\mu_T^*-\mu)\ls z \bigr)\pto\Phi\biggl(
\frac{z}{\sqrt{2\pi f(0)}} \biggr),
\end{equation}
where $\Phi(\cdot)$ denotes the standard normal distribution
function, so that the asymptotic distribution is normal with mean 0 and
variance $2\pi f(0)$.
Then it holds (in probability)
\[
\Biggl\{\frac{1}{\sqrt{2\pi f(0) T}}\sum_{t=1}^{\lfloor T u\rfloor
}\bigl(V^*(t)-\mu\bigr)\dvtx0\ls u\ls1 \Big|  V(\cdot)\Biggr\}\stackrel
{D[0,1]}{\longrightarrow}\{W(u)\dvtx0\ls u\ls1\},
\]
where $V^*(t)=Z^*(t)+\mu_T^*$.
\end{cor}

Along the same lines of the proof we also obtain the analogue of the
finite-sample result of Braun and Kulperger~\cite{braun} for the
surrogate data method. This shows that any finite sample has the same
covariance structure as the corresponding finite sample of the original
sequence; however, it also shows that not only do the partial sums of
the bootstrap sample become more and more Gaussian, but each individual
bootstrap observation does as well.
\begin{cor}\label{cor_braun}
If Assumptions~\ref{ass_W}, \ref
{ass_boot_1}--\ref{ass_boot_4} are fulfilled, then for
any subset $l_1,\ldots,l_p$ of fixed positive integers it holds (in
probability)
\[
(Z^*(l_1),\ldots,Z^*(l_p) |  V(\cdot))\dto
N(0,\Sigma),
\]
where $\Sigma=(\sigma_{i,j})_{i,j=1,\ldots,p}$ with $\sigma
_{i,j}=\cov(V(l_i),V(l_j))$.
\end{cor}


\section{Validity of the meta-assumptions on the bootstrap Fourier
coefficients}
\label{section_valboot}

In this section we prove the validity of the bootstrap schemes if the
Fourier coefficients satisfy certain properties. These or related
properties have been investigated by many researchers in the last
decades and hold true for a large variety of processes. Some of these
results are given in Section~\ref{section_prop_freq}.

Recall Assumption~\ref{ass_a1}, which is important for
the residual-based bootstrap (RB) as well as the wild bootstrap (WB).
\renewcommand{\theassumptionA}{$\mathcal{A}$.\arabic{assumptionA}}
\begin{assumptionA}\label{ass_a1}
Let $\widehat{f}(\cdot)$ estimate the spectral density $f(\cdot)$ in
a uniform way, that is, fulfilling (\ref{eq_est_spec_dens}),
$ {\sup_{\lambda\in[0,\pi]}}|\widehat{f}(\lambda)-f(\lambda
)|\pto0$.
\end{assumptionA}

The next two assumptions are necessary to obtain the validity of the
residual-based bootstrap.
\begin{assumptionA}\label{ass_a2}The following assertions on sums of
the periodogram and/or Fourier coefficients hold:
\begin{eqnarray*}
&\displaystyle \mbox{\textup{(i)}\quad} \frac{1}{2 N}\sum_{j=1}^{N}\frac{x(j)+y(j)}{\sqrt
{f(\lambda_j)}}\pto0, \qquad\mbox{\textup{(ii)}\quad} \frac{1}{2\pi N}\sum
_{j=1}^{N}\frac{I(j)}{f(\lambda_j)}\pto1,&\\
&\displaystyle \mbox{\textup{(iii)}\quad} \frac{1}{ N}\sum_{j=1}^{N}\frac
{I^2(j)}{f^2(\lambda_j)}\ls C+o_P(1)&
\end{eqnarray*}
for some constant $C>0$, where $I(j)=x^{2}(j)+y^2(j)$ is the periodogram.

In particular, \textup{(ii)} is fulfilled if $\frac{1}{ N}\sum
_{j=1}^{N}\frac{I^2(j)}{f^2(\lambda_j)}\pto C>0$.
\end{assumptionA}
\begin{assumptionA}\label{ass_a3}
The empirical distribution function based on the Fourier coefficients
converges uniformly to the standard normal distribution function $\Phi
(\cdot)$, that is,
\[
\sup_{z\in\mathbb{R}}\Biggl|\frac{1}{2N}\sum_{j=1}^N
\bigl(1_{\{x(j)\ls z \sqrt{\pi f(\lambda_j)}\}}+1_{\{
y(j)\ls z \sqrt{\pi f(\lambda_j)}\}}\bigr)-\Phi(z)
\Biggr|\pto0.
\]
\end{assumptionA}

The following two assumptions are necessary to obtain the validity of
the local bootstrap (LB).
\begin{assumptionA}\label{ass_a4}
The following assertions on sums of the periodogram and/or Fourier
coefficients hold true:
\begin{eqnarray*}
\mbox{(i)\hspace*{8pt}}&& \sup_{1\ls k\ls N}\Biggl|\sum_{j=-\infty}^{\infty
}p_{j,T}\bigl(x(k+j)+y(k+j)\bigr)\Biggr|=o_P(1),\\
\mbox{(ii)\hspace*{8pt}}&& \sup_{1\ls k\ls N}\Biggl|\sum_{j=-\infty
}^{\infty}p_{j,T}I(k+j)-2\pi f(\lambda_k)\Biggr|=o_P(1),\\
\mbox{(iii)\hspace*{8pt}}&& \sup_{1\ls k\ls N}\sum_{j=-\infty}^{\infty
}p_{j,T}I^2(k+j)\ls C+o_P(1),
\end{eqnarray*}
where $I(j)=x^2(j)+y^2(j)$ if $j$ is not a multiple of $T$ and
$I(cT)=0$ for $c\in\mathbb{Z}$ and $p_{j,T}$ are as in (\ref
{eq_weights_local}).
\end{assumptionA}
\begin{assumptionA}\label{ass_a5}
The empirical distribution function based on the Fourier coefficients
converges uniformly to the standard normal distribution function $\Phi
(\cdot)$, that is,
\begin{eqnarray*}
&& \sup_s P\Biggl( \sup_{z\in\mathbb{R}}\Biggl|\frac{1}{2}\sum
_{j=1}^Np_{j-s,T}\bigl(1_{\{x(j)\ls z \sqrt{\pi f(\lambda
_j)}\}}\\
&&\hphantom{\sup_s P\Biggl( \sup_{z\in\mathbb{R}}\Biggl|\frac{1}{2}\sum
_{j=1}^Np_{j-s,T}\bigl(}
{}+1_{\{y(j)\ls z \sqrt{\pi f(\lambda_j)}\}
}\bigr)-\Phi(z)\Biggr|\gs\eps\Biggr)\\
&&\qquad \to0,
\end{eqnarray*}
where $p_{j,T}$ are as in (\ref{eq_weights_local}).
\end{assumptionA}

The next theorem shows that the bootstrap methods RB, WB and LB fulfill
the Assumptions~\ref{ass_boot_1}--\ref
{ass_boot_3} which are needed to obtain the functional limit theorem
for the partial sums (cf. Theorem~\ref{th_main}) under the above
assumptions on the periodograms.
\begin{theorem}\label{th_boot}Let Assumption \ref
{ass_W} be fulfilled.
\begin{enumerate}[(a)]
\item[(a)] All three bootstrap methods RB, WB and LB fulfill Assumption~\ref{ass_boot_1} by definition.
\item[(b)] \textit{Residual-based bootstrap RB}: Let Assumption~\ref{ass_process_dens}
hold:
\begin{enumerate}[(iii)]
\item[(i)] Under Assumption~\ref{ass_a1} RB fulfills
Assumption~\ref{ass_boot_2b}.
\item[(ii)] If additionally Assumption~\ref{ass_a2} holds, RB fulfills
Assumption~\ref{ass_boot_4}.
\item[(iii)] If additionally Assumption~\ref{ass_a3} holds,
then Assumption~\ref{ass_boot_3} holds.
\end{enumerate}
\item[(c)] \textit{Wild bootstrap WB}:
under Assumption~\ref{ass_a1} WB fulfills Assumptions \ref
{ass_boot_2b}--\ref{ass_boot_3}.
\item[(d)] \textit{Local bootstrap LB}:
\begin{enumerate}[(iii)]
\item[(i)] Under Assumptions~\ref{ass_a4}\textup{(i)} and~\ref{ass_a4}\textup{(ii)} LB
fulfills Assumption~\ref{ass_boot_2b}.
\item[(ii)] If additionally Assumption~\ref{ass_a4}\textup{(iii)}
holds, then LB fulfills Assumption~\ref{ass_boot_4}.
\item[(iii)] If additionally Assumption~\ref{ass_a5} holds, then LB
fulfills Assumption~\ref{ass_boot_3}.
\end{enumerate}
\end{enumerate}
\end{theorem}

In many applications we apply bootstrap methods not directly to a
stationary sequence $\{V(t)\dvtx t\gs1\}$ but rather to an estimate $\{\wv
(t)\dvtx t\gs1\}$ thereof.
The next corollary gives conditions under which the bootstrap schemes
remain valid in this situation.

In order to give some general conditions under which the procedure
remains valid we need to specify the type of spectral density estimator
we want to use for the residual-based as well as the wild bootstrap. We
want to use the following kernel density estimator (see also Lemma \ref
{lem_spectralestimate}):
%
\begin{equation}\label{eq_lem_spectralestimate}
\widehat{f}_T(\lambda)=\frac{\sum_{j\in\mathbb{Z}}K(
({\lambda-\lambda_j})/{h_T} )I(j)}{2\pi\sum_{j\in\mathbb
{Z}}K({\lambda_j}/{h_T})},
\end{equation}
where $I(j)=x^{2}(j)+y^2(j)$ is the periodogram at frequency $\lambda
_j=2\pi j/T$ if $j$ is not a multiple of $T$ and
$I(cT)=0$, $c\in\mathbb{Z}$.

In addition we need the following assumptions on the kernel.
\renewcommand{\theassumptionK}{$\mathcal{K}$.\arabic{assumptionK}}
\begin{assumptionK} \label{ass_kernelneu1}
Let $K(\cdot)$ be a positive even function with $\int K(\lambda)
\,d\lambda=1$ and
\[
\frac{2\pi}{T h_T}\sum_{j\in\mathbb{Z}}K\biggl( \frac{2\pi j}{T
h_T }\biggr)=\int K(x) \,dx+o(1)=1+o(1).
\]
\end{assumptionK}
\begin{assumptionK}\label{ass_kernelneu2}
Let
\[
{\sup_{\lambda\in[0,2\pi]}}|K_h(\lambda)|=O( h_T^{-1} ),
\]
where
%
\begin{equation} \label{eq_def_kh}
K_h(\lambda)=\frac{1}{h_T}\sum_{j\in\mathbb{Z}}K\biggl( \frac
{\lambda+2\pi j}{h_T} \biggr).
\end{equation}
\end{assumptionK}
\begin{rem}\label{rem_new_kernel}
The above assumption is not as restrictive as it may seem. For example,
it is fulfilled for bounded kernels with compact support. More
generally it holds that (cf., e.g., Priestley~\cite{priestleybook},
equations (6.2.93)--(6.2.95)),
\[
K_h(\lambda)=\frac{1}{2\pi}\sum_{j\in\mathbb{Z}}k( j h_T
)\exp(-ij\lambda),
\]
where $k(x)$ is the inverse Fourier transform of the kernel $K(\cdot
)$, that is,
\[
k(x)=\frac{1}{2\pi}\int K(\lambda)\exp(ix\lambda) \,d\lambda,
\]
respectively,
\[
K(\lambda)=\frac{1}{2\pi}\int k(x)\exp
(-ix\lambda) \,dx.
\]
From the above representation it is clear that as soon as the sum in
$K_h(\lambda)$ can be approximated by an integral for $h_T$ small
enough, it holds for $T$ large
\[
K_h(\lambda)\approx\frac{1}{h_T}K\biggl( \frac{\lambda}{h_T}
\biggr),
\]
which yields the above assumption again for bounded $K(\cdot)$.
Assumption~\ref{ass_kernelneu2} is correct under
Assumptions~\ref{ass_kern_robinson} or~\ref{ass_kern_shaowu} given in the next section.
\end{rem}

We are now ready to state the corollary.
\begin{cor}\label{cor_boot_est}
Assume that the respective (for each bootstrap) conditions of
Theorem~\ref{th_boot} are fulfilled for $\{V(\cdot)\}$.
Moreover assume that we have observed a sequence $\{Y(\cdot)\}$ from
which we can estimate $\{V(\cdot)\}$ by $\{\wv(\cdot)\}$ such that
%
\begin{equation}\label{eq_cond_est}
\frac{1}{T}\sum_{t=1}^T\bigl(V(t)-\wv(t)\bigr)^2=o_P(\alpha_T^{-1})\vadjust{\goodbreak}
\end{equation}
for $\alpha_T\to\infty$ defined below.
Furthermore we assume that for the residual-based and wild bootstrap we
use the spectral density estimator (\ref{eq_lem_spectralestimate})
with a kernel fulfilling Assumptions \ref
{ass_kernelneu1} and~\ref{ass_kernelneu2}. For the
local bootstrap we also use a kernel fulfilling Assumptions \ref
{ass_kernelneu1} and~\ref{ass_kernelneu2}.
\begin{longlist}[(a)]
\item[(a)] If $\alpha_T=T^{1/2}+h_T^{-1}$ [$h_T$ is the bandwidth in (\ref
{eq_lem_spectralestimate})], then the assertions of Theorem \ref
{th_boot} for the residual-based bootstrap RB remain true, but now
given $\{Y(\cdot)\}$.
\item[(b)] If $\alpha_T=h_T^{-1}$ [$h_T$ is the bandwidth in (\ref
{eq_lem_spectralestimate})], then the assertions of Theorem~\ref
{th_boot} for the wild bootstrap WB remain true given $\{Y(\cdot)\}$.

\item[(c)] If $\alpha_T=(T/h_T)^{1/2}$ ($h_T$ is the bandwidth as in the
description of the local bootstrap LB), then the assertions of
Theorem~\ref{th_boot} for the local bootstrap LB remain true given $\{
Y(\cdot)\}$.
\end{longlist}
\end{cor}
\begin{rem}\label{rem_est_err}
There is no assumption on $\{Y(\cdot)\}$ except for (\ref
{eq_cond_est}). In Section~\ref{section_appl} different examples are
given showing the diversity of possible $\{Y(\cdot)\}$ including
nonstationary processes.

Assumption (\ref{eq_cond_est}) is not the weakest possible as the proof
shows. However, it is a condition that is fulfilled in many situations
and one that is easy to verify (cf. Section~\ref{section_appl}).
Weaker conditions would typically include many sine and cosine terms
and would therefore be much more complicated to verify.
\end{rem}

\section{Some properties of Fourier coefficients and periodograms}
\label{section_prop_freq}
In this section we give some examples of processes as well as kernels
which fulfill the assumptions of the previous section. This shows that
the bootstrap has the correct second-order moment structure for a large
class of processes including nonlinear ones.

For the sake of completeness we summarize some recent results of kernel
spectral density estimation leading to Assumption \ref
{ass_a1} in Lemma~\ref{lem_spectralestimate}. In Lemma \ref
{lem_llnperio2}(b) Assumption~\ref{ass_a1} is also proved under a
different set of assumptions. We give now a set of assumptions on the
kernel as well as on the underlying process given by Robinson \cite
{robinson91}, respectively, Shao and Wu~\cite{shaowu07} to obtain
consistent kernel spectral density estimates.
\begin{assumptionK}\label{ass_kern_robinson}
Let the kernel $K(\cdot)$ be a real, even function with
\[
\int_{\mathbb{R}}|K(\lambda)| \,d\lambda<\infty,\qquad
\int_{\mathbb{R}}K(\lambda) \,d\lambda=1.
\]
Furthermore the inverse Fourier transform of $K$
\[
k(x)=\int_{\mathbb{R}}K(\lambda)\exp(ix\lambda) \,d\lambda
\]
satisfies $|k(x)|\ls\widetilde{k}(x)$, where $\widetilde{k}(x)$ is
monotonically decreasing on $[0,\infty)$ and chosen to be an even
function with
\[
\int_0^{\infty}(1+x)\widetilde{k}(x) \,dx<\infty.\vadjust{\goodbreak}
\]
\end{assumptionK}

\begin{assumptionK}\label{ass_kern_shaowu}
Let the kernel $K(\cdot)$ be a real, even function with $\int
_{\mathbb{R}}K(\lambda) \,d\lambda=1$.
Furthermore the inverse Fourier transform of $K$
\[
k(x)=\int_{\mathbb{R}}K(\lambda)\exp(ix\lambda) \,d\lambda
\]
is Lipschitz continuous with support $[-1,1]$.
\end{assumptionK}
\begin{assumptionP}\label{ass_process_robinson}
Assume
\[
\Gamma_V(r)=N^{-1}\sum_{j=1}^{N-r}\bigl(V(j)-\E V(0)\bigr)\bigl(V(j+r)-\E
V(0)\bigr)
\]
satisfies uniformly in $r$ as $T\to\infty$
\[
\Gamma_V(r)-\E\Gamma_V(r)=O_P(N^{-\nu})
\]
for some $0<\nu\ls1/2$.

A detailed discussion of this assumption can be found in Robinson \cite
{robinson91}; for linear processes with existing fourth moments
Assumption~\ref{ass_process_robinson}
is always fulfilled with $\nu=1/2$.
\end{assumptionP}
\begin{assumptionP}\label{ass_process_shaowu1}
Assume that $V(j)-\E V(0)=G(\ldots,\eps(j-1),\eps(j))$, where
$G(\cdot)$ is a measurable function and $\{\eps(\cdot)\}$ is an
i.i.d. sequence. Assume further that
\begin{eqnarray*}
&&\sum_{j=0}^{\infty}\bigl\|\E\bigl(V(j)-\E V(0) | \eps(0),\eps
(-1),\ldots\bigr)- \E\bigl(V(j)-\E V(0) | \eps(-1),\eps(-2),\ldots\bigr)
\bigr\|\\
&&\qquad<\infty,
\end{eqnarray*}
where $\|X\|=\sqrt{\E|X|^2}$.
In case of linear processes this condition is equivalent to the
absolute summability of the coefficients.
\end{assumptionP}

The next assumption is stronger:
\begin{assumptionP}\label{ass_process_shaowu}
Assume that $V(j)=G(\ldots,\eps(j-1),\eps(j))$, where $G(\cdot)$ is
a measurable function and $\{\eps(\cdot)\}$ is an i.i.d. sequence.
Further assume the following geometric-moment contraction condition
holds. Let $\{\widetilde{\eps}(\cdot)\}$ be an i.i.d. copy of $\{\eps
(\cdot)\}$, let $\widetilde{V}(j)=G(\ldots,\widetilde{\eps}(-1),\widetilde
{\eps}(0),\eps(1),\ldots,\eps(j))$ be a coupled version of $V(j)$.
Assume there exist $\alpha>0, C>0$ and $0<\rho=\rho(\alpha)<1$,
such that for all $j\in\mathbb{N}$
\[
\E\bigl( |\widetilde{V}(j)-V(j)|^{\alpha} \bigr)\ls
C\rho^{j}.
\]
This condition is fulfilled for linear processes with finite variance
that are short-range dependent. Furthermore it is fulfilled for a large
class of nonlinear processes. For a detailed discussion of this
condition we refer to Shao and Wu~\cite{shaowu07}, Section~5.
\end{assumptionP}
\begin{assumptionP}\label{ass_process_wu}
Assume that $V(j)=G(\ldots,\eps(j-1),\eps(j))$ is a stationary
causal process, where $G(\cdot)$ is a measurable function and $\{\eps
(\cdot)\}$ is an i.i.d. sequence. Let $\widetilde{V}(j)=G(\ldots,{\eps
}(-1),\widetilde{\eps}(0),\eps(1),\ldots,\eps(j))$ be a coupled
version of $V(j)$ where $\widetilde{\eps}(0)\deq\eps(0)$ independent of
$\{\eps(j)\}$. Furthermore assume
\[
\sum_{i\gs0} \bigl( \E(V_i-\widetilde{V}_i)^2 \bigr)^{1/2}<\infty.
\]
\end{assumptionP}

The following lemma gives some conditions under which
Assumption~\ref{ass_a1} holds, which is necessary for
the residual-based and wild bootstrap (RB and WB) to be valid. Moreover
it yields the validity of Assumption~\ref{ass_a4}(ii),
which is needed for the local bootstrap LB. Lemma~\ref{lem_llnperio2}
also gives some assumptions under which the kernel spectral density
estimate uniformly approximates the spectral density.
\begin{lemma}\label{lem_spectralestimate}
Assume that Assumptions~\ref{ass_W} and~\ref{ass_kernelneu1} are fulfilled, and let $\widehat
{f}_T(\lambda)$ be as in (\ref{eq_lem_spectralestimate}).
\begin{longlist}[(a)]
\item[(a)] Let Assumptions~\ref{ass_kern_robinson} and
\ref{ass_process_robinson} be fulfilled; additionally
the bandwidth needs to fulfill $h_T+h_T^{-1}T^{-\nu}\pto0$. Then
\[
\max_{\lambda\in[0,2\pi]}|\widehat{f}_T(\lambda)-f(\lambda
)|\pto0.
\]
\item[(b)] Let Assumptions~\ref{ass_kern_shaowu},
\ref{ass_clt},~\ref{ass_process_dens} and
\ref{ass_process_shaowu} be fulfilled. Furthermore let
$\E|V(0)|^{\mu}<\infty$ for some $4<\mu\ls8$ and $h_T\to
0$, $(h_T T^{\eta})^{-1}=O(1)$ for some $0<\eta<(\mu-4)/\mu$, then
\[
{\max_{\lambda\in[0,2\pi]}}|\widehat{f}_T(\lambda)-f(\lambda
)|\pto0.
\]
\end{longlist}
\end{lemma}
\begin{rem}
For linear processes Robinson~\cite{robinson91} gives an automatic
bandwidth selection procedure for the above estimator; see also
Politis~\cite{politis03}.
\end{rem}

The following lemma establishes the validity of Assumptions~\ref{ass_a2}.
\begin{lemma}\label{lem_llnperio}Let Assumption \ref
{ass_W} be fulfilled.
\begin{longlist}[(a)]
\item[(a)] Then Assumption~\ref{ass_a2}\textup{(i)} holds.
%
\item[(b)] If additionally
%
\begin{equation}\label{eq_cov_1}
{\max_{1\ls j,k\ls N}}|\cov(I(j),I(k))-(2\pi f(\lambda
_j))^{2}\delta_{j,k}|=o(1),
\end{equation}
where $\delta_{j,k}=1$ if $j=k$ and else $0$, then Assumption
\ref{ass_a2}\textup{(ii)} holds.
%
\item[(c)] If additionally
%
\begin{equation}\label{eq_cov_2}
{\max_{1\ls j,k\ls N}}|\cov(I^{2}(j),I^2(k))-4(2\pi f(\lambda
_j))^{4}\delta_{j,k}|=o(1),\vadjust{\goodbreak}
\end{equation}
then Assumption~\ref{ass_a2}\textup{(iii)} holds, and, more precisely,
\[
\frac{1}{N}\sum_{j=1}^{N}\frac{I^{2}(j)}{f^2(\lambda_j)}=2(2\pi
)^{2}+o_P(1).
\]
\end{longlist}
\end{lemma}
\begin{rem}\label{rem_llnperio}
The conditions of Lemma~\ref{lem_llnperio} are fulfilled for a large
class of processes.
\begin{longlist}
\item Theorem 10.3.2 in Brockwell and Davis~\cite{brockwelldavis}
shows (\ref{eq_cov_1}) for linear processes $V(t)-\E V(t)=\sum
_{j=-\infty}^{\infty}w_j \eps(t-j)$, where $\{\eps(\cdot)\}$ is
i.i.d. with $\E\eps(0)=0$, $\E\eps(0)^4<\infty$ and $\sum
|w_j|\sqrt{j}<\infty$. Furthermore the rate of convergence for $j\neq
k$ is uniformly $O(1/T)$. An analogous proof also yields
(\ref{eq_cov_2}) under the existence of 8th moments, that is, if $\E
\eps(0)^8<\infty$.

\item Lemma A.4 in Shao and Wu~\cite{shaowu07} shows that (\ref
{eq_cov_1}) [resp., (\ref{eq_cov_2})] is fulfilled if the 4th-order
cumulants (resp., 8th-order cumulants) are summable, 4th (resp., 8th)
moments exist and $\sum_{j\gs0}|j \gamma(j)|<\infty$ (cf. also
Theorem 4.3.1 in Brillinger~\cite{brillinger}). More precisely they
show that the convergence rate is $O(1/T)$ in (\ref{eq_cov_1}).
By Remark~4.2 in Shao and Wu~\cite{shaowu07} this cumulant condition
is fulfilled for processes fulfilling Assumption \ref
{ass_process_shaowu} for $\alpha=4$ (resp., $\alpha=8$).
\item Furthermore, Assumption~\ref{ass_a2} is fulfilled
if Assumptions~\ref{ass_process_dens}
and~\ref{ass_process_wu} are fulfilled (cf. Wu~\cite{wu08}).
\item Chiu~\cite{chiu88} uses cumulant conditions to prove strong laws
of large numbers and the corresponding central limit theorems for
$\frac{1}{N}\sum_{j=1}^N\psi(\lambda_j)I^k(\lambda_j)$ for $k\gs1$.
\end{longlist}
\end{rem}

The next lemma shows that weighted and unweighted empirical
distribution functions of Fourier coefficients converge to a normal
distribution, hence showing that Assumptions \ref
{ass_a3} and~\ref{ass_a5} are valid. The proof is based
on Theorem 2.1 in Shao and Wu~\cite{shaowu07}, which is somewhat
stronger than the usual statement on asymptotic normality of finitely
many Fourier coefficients as it gives the assertion uniformly over all
finite sets of fixed cardinal numbers; this is crucial for the proof of
Lemma~\ref{lem_edf}.
\begin{lemma}\label{lem_edf}
Let Assumptions~\ref{ass_W}, \ref
{ass_process_dens} and~\ref{ass_process_shaowu1} be fulfilled.
Furthermore consider weights $\{w_{j,N}\dvtx1\ls j\ls n\}$ such that $\sum
_{j=1}^N w_{j,N}=1$ and\break $\sum_{j=1}^Nw_{j,N}^2\to0$ as $N\to\infty
$, then
\[
\sup_{z\in\mathbb{R}}\Biggl|\frac{1}{2}\sum_{j=1}^Nw_{j,N}
\bigl(1_{\{x(j)\ls z \sqrt{\pi f(\lambda_j)}\}}+1_{\{
y(j)\ls z \sqrt{\pi f(\lambda_j)}\}}\bigr)-\Phi(z)
\Biggr|\pto0,
\]
where $\Phi(\cdot)$ denotes the distribution function of the standard
normal distribution.
If we\vspace*{2pt} have weights $w_{j,N,s}$ with $\sum_{j=1}^N w_{j,N,s}=1$ and
$\sup_s
\sum_{j=1}^Nw_{j,N,s}^2\to0$, then the assertion remains true in the
sense that for any $\eps>0$\vadjust{\goodbreak} it holds that
\begin{eqnarray*}
&&\sup_s P\Biggl( \sup_{z\in\mathbb{R}}\Biggl|\frac{1}{2}\sum
_{j=1}^Nw_{j,N,s}\bigl(1_{\{x(j)\ls z \sqrt{\pi f(\lambda
_j)}\}}+1_{\{y(j)\ls z \sqrt{\pi f(\lambda_j)}\}
}\bigr)-\Phi(z)\Biggr|\gs\eps\Biggr)\\
&&\qquad\to0.
\end{eqnarray*}
\end{lemma}

The next lemma shows the validity of Assumptions~\ref{ass_a4} and again~\ref{ass_a1} under a different set
of assumptions. For this we need to introduce yet another assumption on
the kernel.
\begin{assumptionK}\label{ass_kernelneu3}
Let $K_h(\lambda)$ as in (\ref{eq_def_kh}) fulfill the following
uniform Lipschitz condition ($\lambda_j=2\pi j/T$):
\[
h_T^2| K_h(\lambda_s)-K_h(\lambda_t)|\ls L_K \biggl|\frac
{s-t}{T}\biggr|.
\]
\end{assumptionK}
\begin{rem}
In case of a uniform Lipschitz continuous kernel with compact support,
the assertion is fulfilled for $h_T$ small enough.
For infinite support kernels we still get Assumption \ref
{ass_kernelneu3} as in Remark~\ref{rem_new_kernel} under certain
stronger regularity conditions.
\end{rem}
\begin{lemma}\label{lem_llnperio2}
Let the process $\{V(\cdot)\}$ fulfill Assumptions \ref
{ass_W} and~\ref{ass_process_dens}. Furthermore the
bandwidth fulfills $(h^3 T)^{-1}=o(1)$ and the kernel $K(\cdot)$
fulfills Assumptions~\ref{ass_kernelneu1} and \ref
{ass_kernelneu3} in addition to $1/(T h_T)\sum_j K^2(2\pi j/(T h_T))=O(1)$.
\begin{longlist}[(a)]
\item[(a)] Assumption~\ref{ass_a4}\textup{(i)} holds, if
%
\begin{eqnarray}\label{eq_llnperio_1}
\frac{1}{2 N}\sum_{j=1}^{N}\bigl(|x(j)|+|y(j)|\bigr)&=&O_P(1),\\
\sup_{1\ls l,k\ls N}| \cov(x(l),x(k))-\pi f(\lambda_k)\delta
_{l,k}|&=&O\biggl( \frac{1}{hT} \biggr),\nonumber\\
\label{eq_cov_coef}
\sup_{1\ls l,k\ls N}| \cov(y(l),y(k))-\pi f(\lambda_k)\delta
_{l,k}|&=&O\biggl( \frac{1}{hT} \biggr).
\end{eqnarray}
\item[(b)] Assumption~\ref{ass_a2}\textup{(ii)} together with
(\ref{eq_cov_1}), where the convergence for $j\neq k$ is uniformly of
rate $(h T)^{-1}$, implies Assumption~\ref{ass_a4}\textup{(ii)}
as well as Assumption~\ref{ass_a1} for the spectral density
estimator given in (\ref{eq_lem_spectralestimate}).
\item[(c)] Assumption~\ref{ass_a2}\textup{(iii)} together with
(\ref{eq_cov_2}), where the convergence for $j\neq k$ is uniformly of
rate $(h T)^{-1}$, implies Assumption~\ref{ass_a4}\textup{(iii)}.
\end{longlist}
\end{lemma}
\begin{rem}\label{rem_llnperio2} By the boundedness of the spectral
density (cf. Assumption~\ref{ass_W}) (\ref
{eq_cov_coef}) follows, for example, from Assumption \ref
{ass_a2}(ii).\vadjust{\goodbreak}

If $\sum_{j\gs1}j^{\nu}|\gamma(j)|<\infty$ for
some $\nu>0$ the rate of convergence in (\ref{eq_cov_coef}) is
\[
T^{-\nu}\mbox{ for } 0<\nu<1,\qquad
\frac{\log T}T\mbox{ for }\nu=1, \qquad  T^{-1} \mbox{ for }\nu>1.
\]
For a proof we refer to the supplementary material~\cite{suppA}, proof
of Lemma~\ref{lem_llnperio}. Some conditions leading to (\ref
{eq_cov_1}), respectively, (\ref{eq_cov_2}), with the required
convergence rates can be found in Remark~\ref{rem_llnperio}.
\end{rem}

\section{Some applications}\label{section_appl}
In this section we show that while our procedure still works for the
same class of periodogram-based statistics as the classical frequency
bootstrap methods, we are also able to apply it to statistics that are
completely based on the time domain representation of the observations,
such as the CUSUM statistic for the detection of a change point in the
location model or the least-squares test statistic in unit-root testing.

\subsection{Statistics based on periodograms}\label{sec_stat_perio}

The classical applications of bootstrap methods in the frequency domain
are kernel spectral density estimators (cf. Franke and H\"{a}rdle \cite
{franke}, Paparoditis and Politis~\cite{papapolitis99}) as well as
ratio statistics and Whittle estimators (cf. Dahlhaus and Janas \cite
{dahlhausjanas}, Paparoditis and Politis~\cite{papapolitis99}). This
includes, for example, Yule--Walker estimators for autoregressive processes.

A simple calculation yields
\[
I^{*}(j)=\frac{1}{T}\bigl( x^{*}(j)+y^*(j) \bigr)^2,
\]
where $I^*(j)$ is the periodogram of the TFT-bootstrap time series
$Z^*(\cdot)$ at $\lambda_j$, and $x^*(\cdot)$, $y^*(\cdot)$ are
defined as in Section~\ref{section_descr_boot}.
Comparing that with the original bootstrap procedures for the
periodograms, we realize that for the wild bootstrap we obtain exactly
the same bootstrap periodogram, whereas for the residual-based as well
as local bootstrap we obtain a closely related bootstrap periodogram
but not exactly the same one. The reason is that we did not
simultaneously draw the real and imaginary part of the bootstrap
Fourier coefficient but further exploited the information that real and
imaginary part are asymptotically independent.
Yet, the proofs showing the validity of the bootstrap procedure for the
above mentioned applications go through, noting that the bootstrap's
real part and imaginary part are (conditionally on the data) independent.
The above discussion shows that the procedures discussed in this paper
inherit the advantages as well as disadvantages of the classical
frequency bootstrap procedures.

\subsection{Change-point tests}\label{section_cpa}

In change-point analysis one is interested in detecting structural
changes in time-series such as, for example, a mean change in the
following AMOC (at-most-one-change) location model:
\[
Y(i)=\cases{\mu_1+V(i),&\quad $1\ls i\ls\widetilde{k}$,\cr
\mu_2+V(i),&\quad $\widetilde{k}<i\ls T$,}\vadjust{\goodbreak}
\]
where $V(\cdot)$ is a stationary process with $\E V(0)=0$; $\mu_1$,
$\mu_2$ and $\widetilde{k}$
are unknown.
The question is whether a mean change occurred at some unknown time
$\widetilde{k}$, the so called change-point. This shows that we are
interested in testing
\[
H_0\dvtx \widetilde{k}<T,\qquad \mu_1\neq\mu_2,\qquad  H_1\dvtx \widetilde{k}=T.
\]
Typically, test statistics in this context are based on centered
partial sums such as the well-known CUSUM statistic,
\[
C_T:=\max_{1\ls k\ls T}\Biggl|\frac{1}{\sqrt{T}}\sum
_{j=1}^k\bigl(Y(j)-\bar{Y}_T\bigr)\Biggr|.
\]

\begin{rem}For simplicity we only discuss the classical CUSUM statistic
above. However, extensions to other test statistics in change-point
analysis, such as
\begin{eqnarray*}
C_T^{(1)}(\alpha)&:=&\max_{1\ls k\ls T}\frac{T^{2\alpha
-1/2}}{(k(T-k))^{\alpha}}\Biggl|\sum_{j=1}^k\bigl(Y(j)-\bar{Y}_T\bigr)
\Biggr|,\qquad 0\ls\alpha<\frac1 2,\\
C_T^{(2)}(\beta)&:=&\frac{1}{T}\sum_{k=1}^{T-1}\frac{T^{2\beta
}}{(k(T-k))^{\beta}}\Biggl(\sum_{j=1}^k\bigl(Y(j)-\bar{Y}_T\bigr)
\Biggr)^2,\qquad 0\ls\beta<2,
\end{eqnarray*}
are straightforward using standard techniques of change-point analysis
(cf., e.g., Kirch~\cite{kirchfreq}, proof of Corollary 6.1). This is
not true for extreme-value type test statistics for which stronger
results are needed.
For a detailed discussion of typical test statistics we refer to Cs\"
{o}rg\H{o} and Horv\'{a}th~\cite{csoehor2,csoehor}.
\end{rem}

If $\{V(\cdot)\}$ fulfills Assumption~\ref{ass_clt} we
obtain the following limit under $H_0$ (cf. also Horv\'{a}th \cite
{horlp} and Antoch, Hu{\v{s}}kov{\'{a}} and Pr{\'{a}}{\v{s}}kov{\'{a}}~\cite{anthuslp}):
%
\begin{equation}\label{eq_cpa_new}
\frac{C_T}{\tau}\dto{\sup_{0\ls t\ls1}}|B(t)|,
\end{equation}
where $\{B(\cdot)\}$ is a Brownian bridge and $\tau^2=2\pi f(0)$,
where $f(\cdot)$ is the spectral density of $\{V(\cdot)\}$.

Kirch~\cite{kirchfreq} has already used permutation methods in the
frequency domain to obtain approximations of critical values for
change-point tests. Her idea was to use random permutations of the
Fourier coefficients taking some symmetry properties into account
before back-transforming them to the time domain using the FFT. However,
the covariance structure of a time series is encoded in the variances
of the Fourier coefficients; hence, this structure is destroyed by a
simple permutation.

We will now apply our TFT-bootstrap to obtain critical values for the
above change-point tests. We do not directly observe the process
$V(\cdot)$ since we do not know whether the null hypothesis or the\vadjust{\goodbreak}
alternative holds true; thus, we estimate $Z(\cdot)=V(\cdot)-\bar
{V}_T$ by
\[
\widehat{Z}(t)=Y(t)-\widehat{\mu}_1 1_{\{t\ls\widehat{\widetilde{k}}\}
}-\widehat{\mu}_2 1_{\{t>\widehat{\widetilde{k}}\}},
\]
where (e.g.) $\widehat{\widetilde{k}}=\arg\max\{|{\sum
_{j=1}^k}(Y(j)-\bar{Y}_T)|\dvtx1\ls k\ls T\}$, $\widehat{\mu}
_1=\frac{1}{\widehat{\widetilde{k}}}\sum_{j=1}^{\widehat{\widetilde
{k}}}Y(j)$ and $\widehat{\mu}_2=\frac{1}{T-\widehat{\widetilde
{k}}}\sum_{j=\widehat{\widetilde{k}}+1}^TY(j)$.

The bootstrap statistic is then given by
\[
C_T^*=\Biggl|\frac{1}{\sqrt{T}}\sum_{j=1}^kZ^{*}(j)\Biggr|,
\]
where $\{Z^*(\cdot)\}$ is the TFT-bootstrap sequence defined in
Section~\ref{section_descr_boot}.

The following theorem shows that the (conditional) limit of the
bootstrap statistic is the same as that of the original statistic under
$H_0$ even if the alternative holds true. Hence, the distribution of
the bootstrap statistic is a good
approximation of the null
distribution of the statistic, and the bootstrap critical values are
asymptotically equivalent to the asymptotic critical values under both
the null hypothesis as well as alternatives. This shows that the
asymptotic test and the bootstrap test are asymptotically equivalent.
In the next section a simulation study shows that frequently we get
better results in small samples when using the TFT-bootstrap.
\begin{theorem}\label{th_cpa}
Suppose that the process $V(t)$ fulfills the H\'{a}jek--Renyi
inequality\vspace*{2pt} (cf., e.g., Lemma \textup{B.1} in Kirch~\cite{kirchdiss} for linear
processes),
and let under $H_1$ $(\widehat{\widetilde{k}}-\widetilde{k})/T=O_P({\beta
_T})$, $\beta_T\to0$. Furthermore let the assumptions in Theorem~\ref
{th_boot} hold and
$\alpha_T \max( \beta_T, \log T /T)\to0$ for $\alpha
_T$ as in Corollary~\ref{cor_boot_est}.
Then it holds under $H_0$ as well as $H_1$ for all $x\in\mathbb{R}$
\[
P^*\biggl( \max_{1\ls k\ls T}\frac{C_T^*}{\tau}\ls x \Big|
Y(1),\ldots,Y(T) \biggr)\pto P\Bigl( {\sup_{0\ls t\ls1}}
|B(t)| \ls x\Bigr),
\]
where $\tau$ is as in (\ref{eq_cpa_new}).
This shows that the corresponding bootstrap test (where one calculates
the critical value from the bootstrap distribution) is asymptotically
equivalent to the asymptotic test above.
\end{theorem}
\begin{rem}
The condition $(\widehat{\widetilde{k}}-\widetilde{k})/T=O_P({\beta_T})$ is
fulfilled for a large class of processes with varying convergence rates
$\beta_T$; for certain linear processes we get the best possible rate
$\beta_T=\frac{1}{T}$ (cf., e.g., Antoch, Hu{\v{s}}kov{\'{a}} and Pr{\'{a}}{\v{s}}kov{\'{a}}~\cite{anthuslp}),
but often in the dependent situation the rates are not as good (cf.,
e.g., Kokoszka and Leipus~\cite{kokleip98}).

It is still possible to get the above result under somewhat stronger
assumptions on $\alpha_T$, that is, on the bandwidth $h_T$, if only
weaker versions of the H\'{a}jek--Renyi inequality are fulfilled as,
for example, given in Appendix~B.1 in Kirch~\cite{kirchdiss} for
fairly general processes.\vadjust{\goodbreak}
\end{rem}
\begin{rem}\label{rem_stud}
For practical purposes it is advisable to use some type of studentizing
here. We propose to use the adapted flat-top estimator with automatic
bandwidth choice described\vspace*{1pt} in Politis~\cite{politis05} for the
asymptotic test as well as for the statistic of the original sample.
Let $\widehat{R}(k)=\frac{1}{T}\sum_{t=1}^{T-k}\widehat
{Z}(t)\widehat{Z}(t+k)$,
\[
w(t)=\cases{
1,&\quad $|t|\ls1/2$,\cr
2(1-|t|),&\quad $1/2<|t|<1$,\cr
0,&\quad $|t|\gs1$,}
\]
and the bandwidth $\Lambda_T=2\widehat{\lambda}$, where $\widehat
{\lambda}$ is the smallest positive integer such that
$|\widehat{R}(\widehat{\lambda}+k)/\widehat{R}(0)|<1.4
\sqrt{\log_{10} T/T}$, for $k=1,2,3$.

Then, the estimator is given by
%
\begin{equation}\label{est_flat_top}
\widehat{\tau}{}^2=\max\Biggl(\widehat{R}(0)+2\sum_{k=1}^{\Lambda
_T} w(k/\Lambda_T)\widehat{R}(k),\frac{1}{T(T-1)}\sum_{j=1}^T
\widehat{Z}(t)^2\Biggr).
\end{equation}
The rightmost part in the parenthesis is chosen to ensure positivity
and scale invariance of the estimator.

For a discussion of a related estimator in change-point analysis we
refer to Hu\v{s}kov\'{a} and Kirch~\cite{huskirchconfstud}.

In the bootstrap domain we propose to use an estimator that is closely
related to the bootstrap procedure, namely an estimator based on the
bootstrap periodograms using the same kernel and bandwidth as for the
bootstrap procedure (cf. also (10.4.7) in Brockwell and Davis \cite
{brockwelldavis})
%
\begin{equation}
\label{boot_est_tau}
\widetilde{\tau}^{2*}=p_{0,T} I^*(1)+\sum_{j\gs1} (p_{j,T}+p_{-j,T})I^*(j),
\end{equation}
where $I^*(j)=x^*(j)+y^*(j)$ is the bootstrap periodogram and
\[
p_{s,T}=\frac{K( {2\pi s}/({Th_T}) )}{\sum_j K(
{2\pi j}/({Th_T}) )}.
\]
It can easily be seen using Assumptions \ref
{ass_boot_1}--\ref{ass_boot_4} that
%
\begin{equation}\label{eq_cons_boot_tau}
\E^*(\widetilde{\tau}^{2*}-\tau^2 )^2\pto0,
\end{equation}
if $\sum_j p_{j,T}^2\to0$ which holds under very weak regularity
conditions on the kernel. This shows that the studentized bootstrap
procedure is asymptotically consistent.

This estimator is naturally related to the bootstrap procedure and has
proved to work best in simulations.
This is similar (although maybe for different reasons) to the block
bootstrap for which G\"{o}tze and K\"{u}nsch~\cite{goetzekuensch96}
showed that, in order to obtain second-order correctness of the
procedure, one needs to studentize\vadjust{\goodbreak} the bootstrap statistic with the
true conditional variance of the bootstrap statistic (which is closely
related to the Bartlett estimator), while for the original statistic
one needs to use a different estimator such as the above mentioned
flat-top estimator.
However, a detailed theoretical investigation of which type of
studentizing is best suitable for the TFT-bootstrap is beyond the scope
of this paper.
\end{rem}

\subsection{Unit root testing}\label{section_unitroot}
Unit root testing is a well studied but difficult problem.
Key early references include Phillips~\cite{phillips87}
and Dickey and Fuller~\cite{dickeyfuller79} (see also
the books by Fuller~\cite{fuller96}
and by Hamilton~\cite{hamilton94}).
Nevertheless, the subject is still very much under
investigation (see, e.g., Cavaliere and Taylor \cite
{cavataylor08,cavataylor09a,cavataylor09b}, Chang and Park~\cite
{changpark03}, Park~\cite{park03}, Paparoditis and Politis \cite
{papapolitis05} and the references therein).


The objective here is to test whether a given set of observations
$Y(0),\ldots,\break Y(T)$ belongs to a stationary or a $I(1)$-time series
(integrated of order one), which means that the time series is not
stationary, but its first-order difference $Y(t)-Y(t-1)$ is stationary.
For simplicity we assume that $Y(0)=0$, and we do not consider a
deterministic trend component in this paper. The hypothesis test of
interest can then be stated as
\[
H_0\dvtx \{Y(\cdot)\}\mbox{ is }I(1),\qquad  H_1\dvtx \{Y(\cdot)\} \mbox{ is
stationary}.
\]
Now we note that for
\[
\rho=\lim_{t\to\infty}\frac{\E Y(t)Y(t-1)}{\E Y(t-1)^2}
\]
the null hypothesis is equivalent to $\rho=1$ (for a detailed
discussion we refer to Paparoditis and Politis~\cite{papapolitis03},
Example 2.1).
Denote $V(t)=Y(t)-\rho Y(t-1)$, which is a stationary sequence under
$H_0$ as well as $H_1$. While the bootstrap test below is valid for the
general situation ($H_0,H_1$) above, it is intuitively easier to
understand if one considers the following restricted situation, where
$Y(t)=\rho Y(t-1)+V(t)$ for some stationary $\{V(\cdot)\}$ with mean 0
and tests $H_0\dvtx \rho=1$ versus $H_1\dvtx|\rho|<1$; this is the setup
we use in the simulations below.

An intuitive test statistic (cf. Phillips~\cite{phillips87}) is given
by $U_T:=T(\widehat{\rho}_T-1)$ rejecting the null hypothesis if
$U_T<c_{\alpha}$ for some appropriate critical value $c_{\alpha}$, where
\[
\widehat{\rho}_T=\frac{\sum_{t=1}^TY(t)Y(t-1)}{\sum_{t=1}^TY(t-1)^2}
\]
is a consistent estimator for $\rho$ under both the null hypothesis as
well as the alternative.
Other choices for $\rho$ and $\widehat{\rho}_T$ are also possible
(for a detailed discussion, see Section 2 of Paparoditis and
Politis~\cite{papapolitis03}).

If $\{V(\cdot)\}$ fulfills Assumption~\ref{ass_clt} with mean 0 and additionally
\[
\frac1 T \sum_{j=1}^T V^2(j)\to
\sigma^2,\vadjust{\goodbreak}
\]
then it holds under $H_0$ that
\[
U_T:=T(\widehat{\rho}_T-1) \dto\frac{W(1)^2-\sigma^2/\tau^2}{2
\int_0^1W(t)^2 \,dt},
\]
%
where $\tau^2=2\pi f(0)$ and $f(\cdot)$ is the spectral density of
the stationary sequence $Y(t)-Y(t-1)$, $\sigma^2=\var(V(t))$, and $\{
W(\cdot)\}$ is a Wiener process (see also Phillips~\cite{phillips87},
Theorem 3.1).

This shows that the limit distribution of $U_T$ depends on the unknown
parameters $\sigma^2$ as well as $\tau^2$ if the errors are
dependent. The famous Dickey--Fuller test is closely related to the
above test statistic just using a slightly different normalization, but
it suffers from the same problem.

Phillips~\cite{phillips87} and Phillips and Perron \cite
{phillipsperron88} suggest some modifications of the two tests
mentioned above which do have a pivotal limit for time series errors as well.
Later on, Perron and Ng~\cite{perronng96}, Stock~\cite{stock99} and
Ng and Perron~\cite{ngperron01} propose to use the trinity of
so-called $M$ unit root statistics, which are also closely related to
the above two tests but have pivotal limits for time series errors as
well. Those $M$ unit root tests are given by
\begin{eqnarray*}
M U_T&=&\frac{(1/T) Y(T)^2-\widehat{\tau}^2}{2/{T^{2}} \sum
_{t=1}^TY(t-1)^2},\\
MSB_T&=&\Biggl( \frac{1}{\widehat{\tau}^2}\frac1 {T^2} \sum
_{t=1}^T Y(t-1)^2\Biggr)^{1/2}
\end{eqnarray*}
as well as the product of the above two $M$ statistics. As before
$\widehat{\tau}^2$ denotes an estimator of $\tau^2$. In the
simulations we use the estimator as given in (\ref{est_flat_top}).

All of the above mentioned statistics are continuous functions of the
partial sum process $\{\sum_{i=1}^{\lfloor T t\rfloor}V(t)\dvtx 0\ls t\ls
1\}$ under the null hypothesis [as $Y(t)=\sum_{j=1}^t V(j)$], so that
the null asymptotics are immediate consequences of the functional
central limit theorem as given in Assumption \ref
{ass_clt}. For the test statistic $U_T$ and the Dickey--Fuller test it
is additionally needed that $\frac1 T \sum_{j=1}^TV^2(j)\to\sigma
^2$. For example, the statistic $M U_T$ has the same asymptotic limit
as the statistic $U_T$ with independent errors
%
\begin{equation}\label{eq_asym_MUT}
M U_T \dto\frac{W(1)^2-1}{2 \int_0^1W(t)^2 \,dt}.
\end{equation}
In the following we concentrate on the statistic $M U_T$ but the
results for the other mentioned statistics follow analogously.

We would like to apply the TFT-bootstrap to obtain critical values;
that means we need a bootstrap sequence which is (conditionally)
$I(1)$. In order to obtain this we estimate $V(t)=Y(t)-\rho Y(t-1)$,
which is stationary under both $H_0$ as well as $H_1$, by
\[
\widehat{V}(t)=Y(t)-\widehat{\rho}_T Y(t-1).
\]
Then we can use the TFT-bootstrap based on $\widehat{V}(\cdot)$, that
is, create
a TFT-bootstrap sample $V^*(1),V^*(2),\ldots$ and obtain a bootstrap
$I(1)$ sequence (i.e., a sequence fulfilling $H_0$) by letting
\[
Y^*(0)=Y(0), \qquad Y^*(t)=Y^*(t-1)+V^*(t), \qquad t\gs1.
\]
The bootstrap analogue of the statistic $MU_T$ is then given by
%
\[
MU_T^*=\frac{(1/T) Y^*(T)-\widetilde{\tau}^{2*}}{2/{T^2}
\sum{t=1}^T(Y^*(t-1))^2},
\]
where we use again the estimator $\widetilde{\tau}^{2*}$ as in (\ref
{boot_est_tau}) for the bootstrap sequence.

The following theorem shows that the conditional limit of the bootstrap
statistic is the same as
that appearing in the RHS of (\ref{eq_asym_MUT}) no matter whether the
original sequence follows the null or alternative hypothesis. This
shows that the bootstrap critical values and thus also the bootstrap
test is equivalent to the asymptotic critical values (and thus the
asymptotic test), under both the null hypothesis as well as alternatives.
\begin{theorem}\label{th_unitroot}
Suppose that the process $\{V(\cdot)\}$ has mean 0 and fulfills the
assumptions in Theorem~\ref{th_boot} and let (\ref{eq_boot_mean}) be
fulfilled. Furthermore assume that under $H_0$ as well as $H_1$ it
holds that
%
\begin{equation}\label{eq_unit_est}
(\rho-\widehat{\rho}_T)^2 \frac1 T \sum
_{t=0}^{T-1}Y^2(t)=o_P(\alpha_T^{-1})
\end{equation}
for $\alpha_T$ as in Corollary~\ref{cor_boot_est}.
Then it holds under $H_0$ as well as $H_1$ for all $x\in\mathbb{R}$ that
%
\[
P^*( MU_T^*\ls x )\pto P\biggl( \frac{W(1)^2-1}{2 \int
_0^1W(t)^2 \,dt} \ls x\biggr).
\]
This shows that the corresponding bootstrap test (where one calculates
the critical value from the bootstrap distribution) is asymptotically
equivalent to the test based on (\ref{eq_asym_MUT}).
\end{theorem}
\begin{rem}
Condition (\ref{eq_unit_est}) is fulfilled for a large class of
processes and if $\alpha_T/T\to0$, Theorem 3.1 in Phillips \cite
{phillips87}, for example, shows under rather general assumptions that
under $H_0$
\[
\widehat{\rho}_T-\rho=O_P(T^{-1}),\qquad \frac{1}{T^2}\sum
_{j=0}^{T-1}Y^2(t)=O_P(1).
\]
Under $H_1$ (\ref{eq_unit_est}) also holds under fairly general
assumptions (cf., e.g., Romano and Thombs~\cite{romanothombs96},
Theorem 3.1) if $\alpha_T/T\to0$; more precisely
\[
\widehat{\rho}_T-\rho=O_P(T^{-1/2}), \qquad\frac{1}{T}\sum
_{j=0}^{T-1}Y^2(t)=O_P(1).
\]
%
\end{rem}

\section{Simulations}\label{section_sim}
In the previous sections, the asymptotic applicability of the
TFT-bootstrap was investigated. In this section we conduct a small
simulation study in order to show its applicability in finite samples.

%
\begin{figure}[t!]

\includegraphics{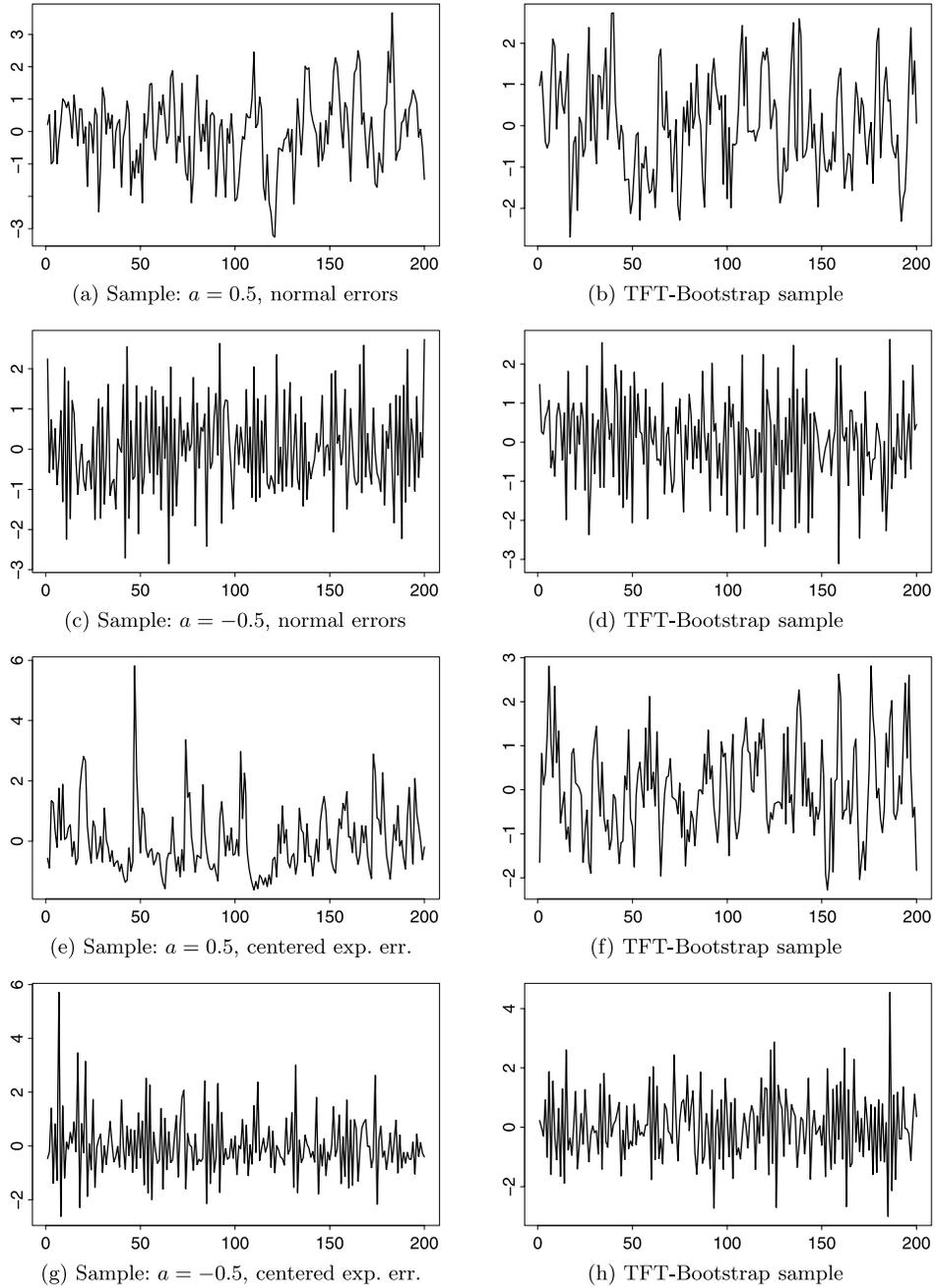}

\caption{$\mathrm{AR}(1)$ with parameter $a$ and corresponding
TFT-bootstrap sample: residual-based bootstrap, Bartlett--Priestley
kernel, $h=0.01$.} \label{fig_1}
\end{figure}

To get a first impression of what a TFT-bootstrap sample looks like, we
refer to Figure~\ref{fig_1}, which shows the original time series as
well as one bootstrap sample. At first glance the covariance structure
is well preserved.

We use the statistical applications given in Section \ref
{section_appl} to show that the TFT-bootstrap is indeed applicable.
The usefulness of the procedure for statistics based on periodograms
have already been shown by several authors (cf. Franke and H\"
{a}rdle~\cite{franke}, Dahlhaus and Janas~\cite{dahlhausjanas} and
Paparoditis and Politis~\cite{papapolitis99}) and will not be
considered again.

However, the applicability for statistics that are completely based on
time domain properties, such as the CUSUM statistic in change-point
analysis or the above unit-root test statistics, is of special
interest. More precisely we will compare the size and power of the
tests with different parameters as well as give a comparison between
the TFT, an asymptotic test, and alternative block resampling
techniques. For change-point tests, the comparison is with the block
permutation test of Kirch~\cite{kirchblock}; in the unit-root
situation we compare
the TFT-bootstrap to the block bootstrap of Paparoditis and
Politis~\cite{papapolitis03}. For the TFT we use the local bootstrap
(LB) as well as residual-based bootstrap (RB) with a uniform kernel
(UK) as well as Bartlett--Priestley kernel (BPK) with different bandwidths.

We visualize these qualities by the following plot:

\subsection*{Achieved size-power curves (ASP)}

The line corresponding to the null hypothesis shows the actual achieved
level on the $y$-axis for a nominal one as given by the $x$-axis. This can
easily be done by plotting the empirical distribution function (EDF) of
the $p$-values of the statistic under $H_0$. The line corresponding to
the alternative shows the size-corrected power, that is, the power of
the test belonging to a true level $\alpha$ test where $\alpha$ is
given by the $x$-axis. This can easily be done by plotting the EDF of
the $p$-values under the null hypothesis against the EDF of the
$p$-values under the alternative.\looseness=1

In the simulations we calculate all bootstrap critical values based on
1,000 bootstrap samples, and the ASPs are calculated on the basis of
1,000 repetitions.

Concerning the parameters for the TFT-bootstrap we have used a uniform
[$K(t)=\frac1 2 1_{\{|t|\ls1\}}$] as well as Bartlett--Priestley
kernel [$K(t)=\frac3 4 (1-t^2)1_{\{|t|\ls1\}}$] with various
bandwidths. All errors are centered exponential hence non-Gaussian.
Furthermore $\mathrm{AR}(1)$ time series are used with coefficient $\gamma
=-0.5,0.5$. Furthermore we consider GARCH$(1,1)$ processes as an example
of nonlinear error sequences.

\subsection{Change-point tests}

\begin{figure}

\includegraphics{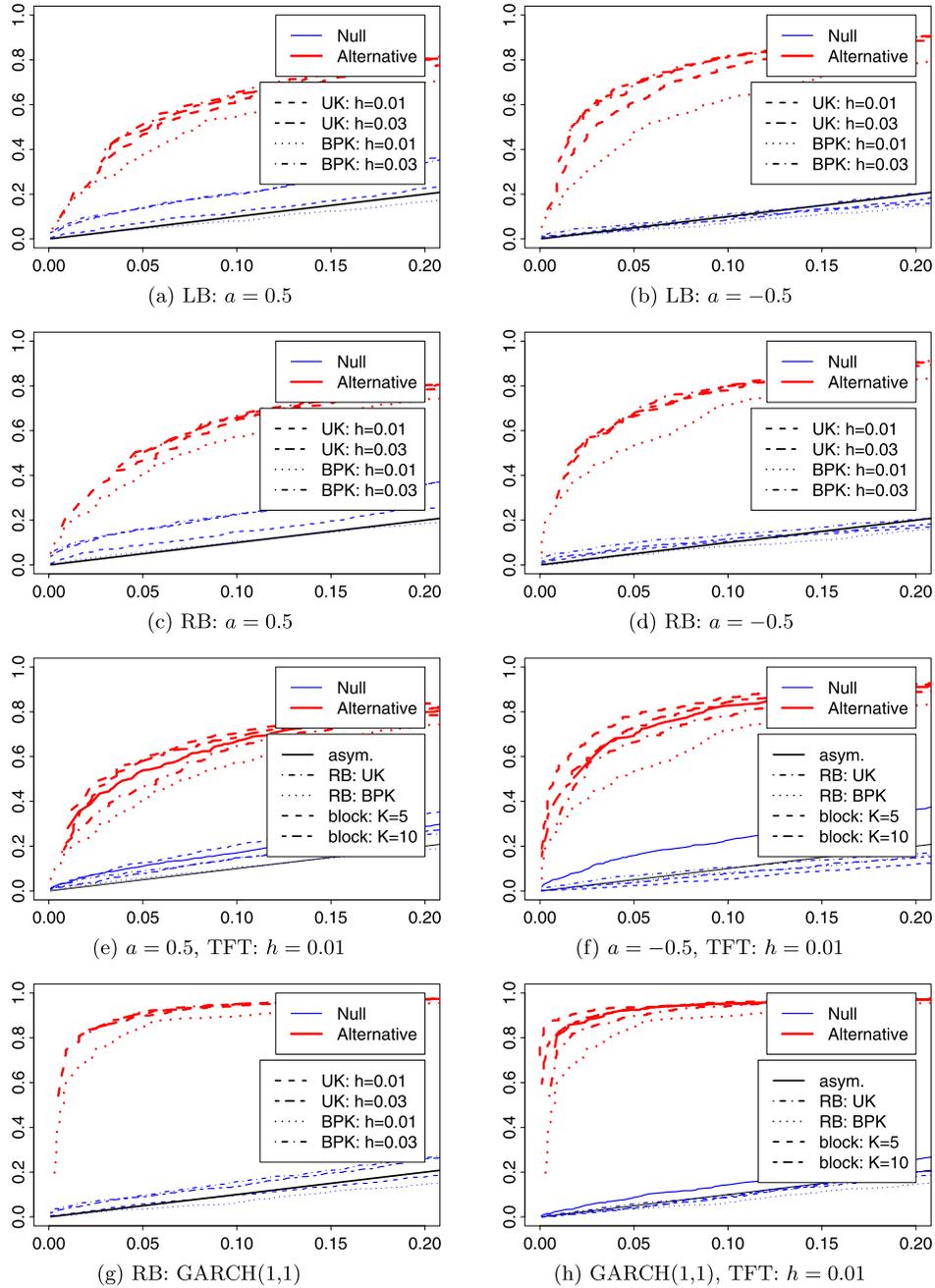}

\caption{Change-point tests: ASP-plots for $\mathrm{AR}(1)$ process
with parameter $a$ with centered exponential errors, respectively,
$\mathrm{GARCH}(1,1)$, $d=\mu_2-\mu_1$, $T=200$, $\widetilde{k}=100$, BPK:
Bartlett--Priestley kernel, UK: uniform kernel, LB: local bootstrap, RB:
residual-based bootstrap.} \label{fig_2}
\end{figure}

We compare the power using an alternative that is detectable but has
not power one already in order to pick up power differences. For the
$\mathrm{AR}(1)$ process with parameter $a=-0.5$, we choose $d=0.3$; for $a=0.5$
we choose $d=0.7$ as changes are more difficult to detect for these
time series. A comparison involving the uniform kernel (UK), the
Bartlett--Priestley kernel (BPK) as well as bandwidth $h=0.01$ and
$h=0.03$ can be found in Figure \mbox{\ref{fig_2}(a)--\ref{fig_2}(d)}. It becomes clear
that small bandwidths are best in terms of keeping the correct size,
where the BPK works even better than the UK. However, this goes along
with a loss in power which is especially severe for the BPK.
Furthermore, the power loss for the BPK kernel is worse if combined
with the local bootstrap. Generally speaking, the TFT works better for
negatively correlated errors which is probably due to the fact that the
correlation between Fourier coefficients is smaller in that case.

%

In a second step, we compare the residual-based bootstrap (RB) with
both kernels and bandwidth $h=0.01$ with the block permutation method
of Kirch~\cite{kirchblock} as well as the asymptotic test. The results
are given in Figure~\ref{fig_2}(e)--\ref{fig_2}(f). The TFT works best in terms
of obtaining correct size, where the BPK beats the UK as already
pointed out above. The power loss of the BPK is also present in
comparison with the asymptotic as well as the block permutation
methods; the power of the uniform kernel is also smaller than for the
other method but not as severe as for the BPK. The reason probably is
the sensitivity of the TFT with respect to the estimation of the
underlying stationary sequence as in Corollary~\ref{cor_main}. In this
example a mis-estimation of the change-point or the mean difference can
result in an estimated sequence that largely deviates from a stationary
sequence, while in the unit-root example below, this is not as
important and in fact the power loss does not occur there.

The\vspace*{-1pt} simulation results for a GARCH$(1,1)$ time series with parameters
$\omega=0.3,\alpha=0.7, \beta=0.2$ [i.e., $V(t)=\sigma_t\eps_t$,
$\sigma_t^2=\omega+\alpha\eps_{t-1}^2+\beta\sigma_{t-1}^2$, $\eps
_t\stackrel{\mathrm{i.i.d.}}{\sim} N(0,1)$] are given in Figure
\ref{fig_2}(g) and~\ref{fig_2}(h), and it becomes clear that the conclusions are similar. The
alternative in these plots is given by $d=0.7$.

\subsection{Unit root testing}

In the unit root situation we need a bootstrap sample of $V(\cdot)$ as
in Corollary~\ref{cor_main}, where we additionally use the fact that
$\{V(\cdot)\}$ has mean $0$. In this case we additionally need a
bootstrap version of the mean. For simplicity we use a wild bootstrap
$\mu^*=W (2\pi\widehat{\tau}^2/T)^{1/2} $, where $W$ is standard
normal distributed and $\widehat{\tau}^2$ is as in (\ref
{est_flat_top}), where we replace $\widehat{Z}(\cdot)$ by $\widehat
{V}(\cdot)-\bar{\widehat{V}}_T$. The alternative in all plots is
given by $\rho=0.95$.

%
\begin{figure}

\includegraphics{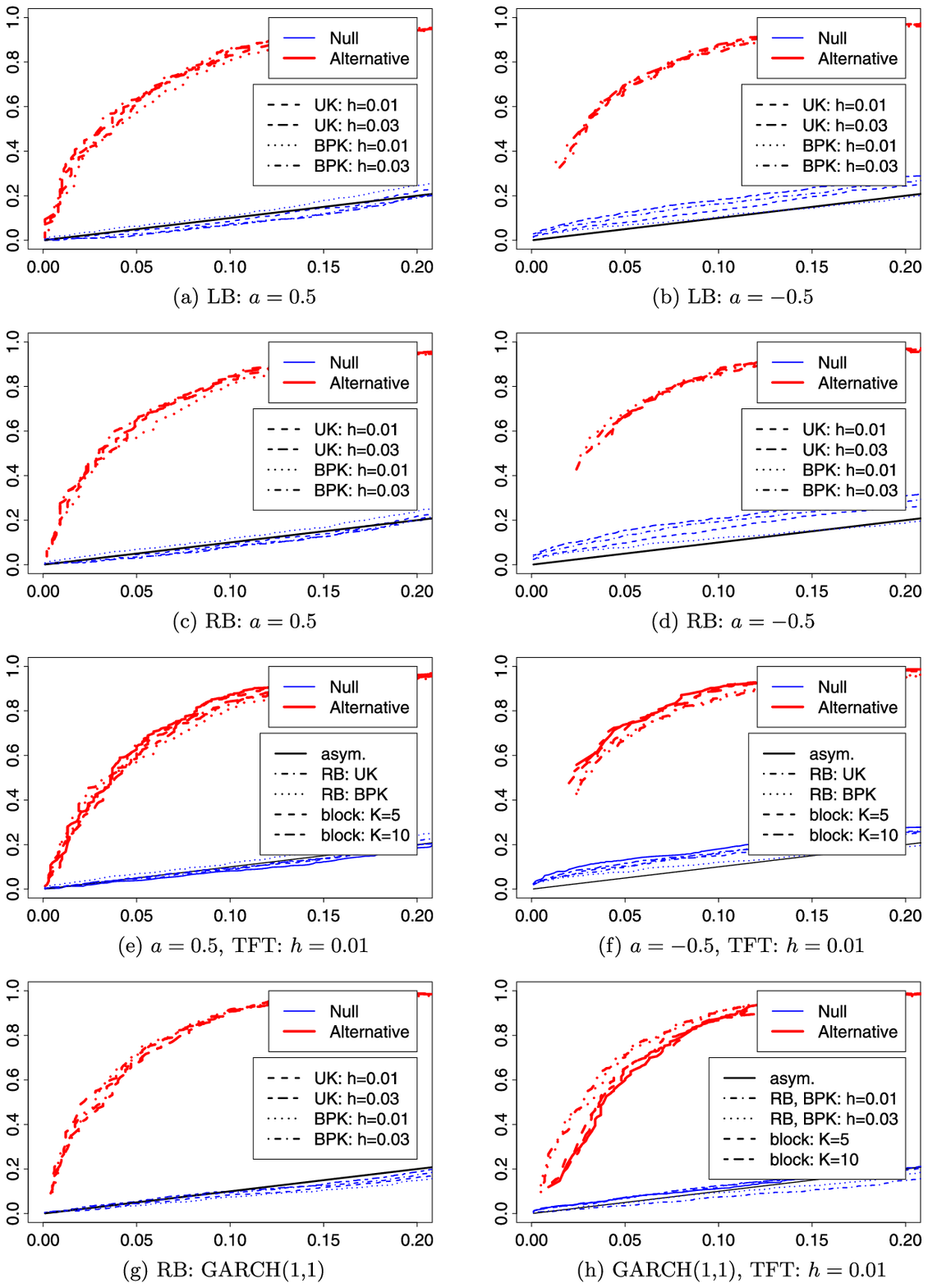}

\caption{Unit-root\vspace*{1pt} tests: ASP-plots for $\mathrm{AR}(1)$ process with
parameter $a$ with centered exponential errors, respectively,
$\mathrm{GARCH}(1,1)$, $d=\mu_2-\mu_1$, $T=200$, $\widetilde{k}=100$, BPK:
Bartlett--Priestley kernel, UK: uniform kernel, LB: local bootstrap,
RB: residual-based bootstrap.} \label{fig_3}
\end{figure}

In Figure~\ref{fig_3}(a)--\ref{fig_3}(d) a comparison of different kernels and
bandwidths for an $\mathrm{AR}(1)$ error sequence is given. It can be seen that
again a small bandwidth yields best results, and in particular the BPK
works better than the UK. Unlike in the change-point example we do not
have the effect of a power loss. Furthermore unlike in change-point
analysis the bootstrap works better for positively correlated errors.

A comparison with the asymptotic test as well as the block bootstrap by
Paparoditis and Politis~\cite{papapolitis03} can be found in
Figure~\ref{fig_3}(e)--\ref{fig_3}(f). In the case of a positive correlation all
methods perform approximately equivalently, at least if we use the
better working Bartlett--Priestley kernel; however for a negative
correlation the TFT holds the level better than the other methods.

Some results for a GARCH$(1,1)$ error sequence with parameters $\omega
=0.3,\alpha=0.7, \beta=0.2$ are shown in Figure~\ref{fig_3}(g) and~\ref{fig_3}(h).
In this situation a somewhat larger bandwidth of $h=0.03$ works
slightly better and the TFT test leads to an improvement of the power
of the tests.

It is noteworthy that the appropriate bandwidth in all cases is smaller
than what one might have expected to be a good choice. A possible
explanation for this is that some undersmoothing is appropriate since
the back-transformation will account for some additional smoothing.

\section{Conclusions}\label{section_conclusions}
The subject of the paper is the TFT-bootstrap which is a
general frequency domain bootstrap method
that also generates bootstrap data in the time domain.
Connections of the TFT-bootstrap
with other methods including the surrogate data method of Theiler et
al.~\cite{theileretal92},
and the original proposals of Hurvich and Zeger~\cite{hurvichzeger}
were thoroughly explored.

It was shown that the TFT-bootstrap samples have asymptotically the
same second-order moment structure as the original time series.
However, the
bootstrap pseudo-series are (asymptotically) Gaussian showing that the
TFT-bootstrap approximates the original time series by a Gaussian
process with the same covariance structure even when the original data
sequence is nonlinear; see Section~\ref{section_fclt}.

Nevertheless, our simulations suggest that for small samples
the TFT-bootstrap gives a better approximation to the critical values
(as compared with the asymptotic ones) especially when studentizing is possible.
Whether appropriate studentization results in higher-order
correctness is a subject for future
theoretical investigation. Choosing somewhat different types of
bootstrapping in the frequency domain could also lead to higher-order
correctness without bootstrapping as in, for example, Dahlhaus and
Janas~\cite{dahlhausjanas2} for the periodogram bootstrap.

In fact the simulations suggest that a surprisingly small bandwidth (in
comparison to spectral density estimation procedures) works best. When
applied in a careful manner no smoothing at all still results---in
theory---in a correct second-moment behavior in the time domain (cf.
Section~\ref{section_diff_methods}), suggesting that due to the
smoothing obtained by the back-transformation a coarser approximation
in the frequency domain is necessary to avoid oversmoothing in the time domain.

\section{Proofs}\label{section_fclt_proof}
In this section we only give a short outline of some of the proofs. All
technical details can be obtained as electronic supplementary
material~\cite{suppA}.

We introduce the notation $a_n\preceq b_n \dvtx\Leftrightarrow a_n=O(b_n)$.

The following lemma is needed to prove Lemma~\ref{lem_cov}.
\begin{lemma}\label{lem_proof_1}
Under Assumption~\ref{ass_W}, the following
representation holds:
\[
\frac{4\pi}{T}\sum_{k=1}^{N}f( \lambda_k )\cos
(\lambda_k h)=\sum_{l\in\mathbb{Z}}\gamma(h+lT)-\frac{2\pi
}{T}f(0) -\frac{2\pi}{T}f(\pi)\exp(i\pi h)1_{\{T\ \mathrm{even}\}}.
\]
\end{lemma}
\begin{pf*}{Proof of Lemma~\ref{lem_cov}}
By Lemma A.4 in Kirch~\cite{kirchfreq} it holds (uniformly in~$u$) that
\[
\sum_{l=1}^{\lfloor m u\rfloor}\cos(\lambda_k l)=O\biggl( \min
\biggl( \frac{T}{k},m\biggr) \biggr).
\]
Thus it holds uniformly in $u$ and $v$ that
%
\begin{equation}\label{eq_sum_trig}
\sum_{k=1}^N\Biggl|\sum_{l_1=1}^{\lfloor m u\rfloor}\cos(\lambda
_kl_1)\Biggr|\Biggl|\sum_{l_2=1}^{\lfloor m v\rfloor}\cos(\lambda
_kl_2)\Biggr|=O(1)\sum_{k=1}^{N}\min\biggl(\frac T k, m
\biggr)^2=O(m T).\hspace*{-32pt}
\end{equation}
By Assumptions~\ref{ass_boot_1} and \ref
{ass_boot_2b} and by (\ref{eq_form_bs})
it holds that
\begin{eqnarray*}
&&\cov^*\Biggl( \sum_{l_1=1}^{\lfloor mu\rfloor}Z^*(l_1),\sum
_{l_2=1}^{\lfloor mv\rfloor}Z^*(l_2)\Biggr)\\
&&\qquad =\frac{4\pi}{T}\sum_{l_1=1}^{\lfloor mu\rfloor}\sum
_{l_2=1}^{\lfloor mv\rfloor}\sum_{k=1}^{N}f(\lambda_k)\cos\bigl(\lambda
_k(l_1-l_2)\bigr)+o_P(m),
\end{eqnarray*}
where the last line follows for $m/T\to0$ as well as $m=T$ by (\ref
{eq_sum_trig}). Assertion (b) follows by an application of Lemma \ref
{lem_proof_1} as well as standard representation of the spectral
density as sums of auto-covariances (cf., e.g., Corollary 4.3.2 in
Brockwell and Davis~\cite{brockwelldavis}).
\end{pf*}

The next lemma gives the crucial step toward tightness of the partial
sum process.
\begin{lemma}\label{lem_tight}
Under Assumptions~\ref{ass_W}, \ref
{ass_boot_1}--\ref{ass_boot_4} it holds for $u<v$ and
some constant $D>0$ that
\[
\E^*\Biggl( \frac{1}{\sqrt{m}}\sum_{l=\lfloor m u\rfloor
+1}^{\lfloor m v\rfloor}Z^*(l) \Biggr)^4\ls\bigl(D+o_P(1)\bigr)(v-u)^2.
\]
\end{lemma}
\begin{pf}
Equation (\ref{eq_sum_trig}) and Lemma~\ref{lem_proof_1} yield the assertion
after some standard calculations.
\end{pf}

The next lemma gives the convergence of the finite-dimensional distribution.
\begin{lemma}\label{lem_finite}
Let $S^*_m(u)=\frac{1}{\sqrt{m}}\sum_{j=1}^{\lfloor m u\rfloor}Z^*(j)$.
\begin{longlist}[(a)]
\item[(a)] If Assumptions~\ref{ass_W}, \ref
{ass_boot_1}--\ref{ass_boot_4} are fulfilled and
$m/T\to0$ we obtain for all $0<u_1,\ldots,u_p\ls1$ in probability
\[
(S_m^*(u_1),\ldots,S_m^*(u_p))\dto N(0,\Sigma),
\]
where $\Sigma=(c_{i,j})_{i,j=1,\ldots,p}$ with $c_{i,j}=2\pi f(0)\min
(u_i,u_j)$.
\item[(b)] If Assumptions~\ref{ass_W}, \ref
{ass_boot_1} and~\ref{ass_boot_3} are fulfilled we obtain
for all $0<u_1,\ldots,u_p\ls1$ in probability
\[
(S_T^*(u_1),\ldots,S_T^*(u_p))\dto N(0,\Sigma),
\]
where $\Sigma=(c_{i,j})_{i,j=1,\ldots,p}$ with $c_{i,j}=2\pi
f(0)(\min(u_i,u_j)-u_i u_j)$.
\end{longlist}
\end{lemma}
\begin{pf}
For assertion (a) we use the Cram\'{e}r--Wold device and prove a
Lyapunov-type condition. Again arguments similar to (\ref
{eq_sum_trig}) are needed.
To use this kind of argument it is essential that $m/T\to0$ because
for $m=T$ the Feller condition is not fulfilled, and thus the Lindeberg
condition can also not be fulfilled.
Therefore a different argument is needed to obtain asymptotic normality
for $m=T$. We make use of the Cram\'{e}r--Wold device and Lemma~3 in
Mallows~\cite{mallows72}. As a result somewhat stronger assumptions
are needed, but it is not clear whether they are really necessary (cf.
also Remark~\ref{rem_lind}).
\end{pf}
\begin{pf*}{Proof of Theorem~\ref{th_main}}
Lemmas~\ref{lem_tight} and~\ref{lem_finite} ensure convergence of the
finite-dimensional distribution as well as tightness, which imply by
Billingsley~\cite{billv2}, Theorem 13.5,
\begin{eqnarray*}
&& \Biggl\{\frac{1}{\sqrt{m}}\sum_{l=1}^{\lfloor mu\rfloor
}\bigl(Z^*(l)-\E^*Z^*(l)\bigr)\dvtx0\ls u\ls1\Biggr\}\\
&&\qquad\stackrel{D[0,1]}{\longrightarrow}
\cases{
\{W(u)\dvtx0\ls u\ls1\}, &\quad $\dfrac m T\to0$,\vspace*{2pt}\cr
\{B(u)\dvtx0\ls u\ls1\},&\quad $m=T$.}
\end{eqnarray*}
\upqed\end{pf*}
\begin{pf*}{Proof of Theorem~\ref{th_boot}}
Concerning Assumption~\ref{ass_boot_2b} it holds by
Assumption~\ref{ass_a1} that
\[
{\sup_k}|{\var}^*(x^*(k))-\pi f(\lambda_k)|={\sup_k}
|\pi\widehat{f}(\lambda_k)-\pi f(\lambda_k)|\pto0.
\]
Similarly we obtain Assumption~\ref{ass_boot_4} since
\[
\sup_k\E^*(x^*(k))^4=3 \pi^2\sup_k\widehat{f}(\lambda_k)^2\ls
3\pi^2\sup_k f(\lambda_k)^2+o_P(1)\ls C+o_P(1).
\]
Concerning Assumption~\ref{ass_boot_3} let $X\deq N(0,1)$, then
$\sqrt{\pi\widehat{f}(\lambda_k)}X\ddeq x^*(k)$. Then
\begin{eqnarray*}
\sup_k d_2^2(\mathcal{L}^*(x^*(k)),N(0,\pi f(\lambda_k)))&\ls&\pi
\sup_k \bigl(\sqrt{f(\lambda_k)}-\sqrt{\widehat{f}(\lambda
_k)}\bigr)^2 \E X^2\\
&\ls&{\pi\sup_k} {|\widehat{f}( \lambda_k)-f(\lambda_k)
|}=o_P(1).
\end{eqnarray*}
The proofs especially of the last statement for the residual-based as
well as local bootstrap are technically more complicated but similar.
\end{pf*}
\begin{pf*}{Proof of Corollary~\ref{cor_boot_est}}
We put an index $V$, respectively, $\wv$, on our previous notation
indicating whether we use $V$ or $\wv$ in the calculation of it, for example,
$x_{\wv}(j)$, $x_V(j)$, respectively, $y_{\wv}(j)$, $y_V(j)$ denote
the Fourier coefficients based on $\wv(\cdot)$, respectively, $V(\cdot
)$.

We obtain the assertion by verifying that Assumptions~\ref{ass_a1},~\ref{ass_a2} as well as
Assumption~\ref{ass_a4} remain true. This in turn implies Assumptions
\ref{ass_boot_2b} as well as Assumption \ref
{ass_boot_4}. Concerning Assumption~\ref{ass_boot_3} we show that
the Mallows distance between the bootstrap r.v. based on $\wv(\cdot)$
and the bootstrap r.v. based on $V(\cdot)$ converges to 0.

The key to the proof is
%
\begin{eqnarray}\label{eq_pcor_boot_3}
&&I_V(j)-I_{\wv}(j)\nonumber\\
&&\qquad=-\frac{1}{ T}\bigl(F^2_T(j) +F^2_T(N+j)\bigr)
+ 2\frac{1}{\sqrt{T}} x_V(j)F_T(j)\\
&&\qquad\quad{} + 2\frac{1}{\sqrt{T}} y_V(j)F_T(N+j),\nonumber
\end{eqnarray}
where
\[
F_T(j):=
\cases{\displaystyle \sum_{t=1}^T\bigl(V(t)-\wv(t)\bigr)\cos(t\lambda_j), &\quad $1\ls j\ls
N$,\vspace*{2pt}\cr
\displaystyle \sum_{t=1}^T\bigl(V(t)-\wv(t)\bigr)\sin(t\lambda_{j-N}), &\quad $N<j\ls2N$.}
\]
This can be seen as follows:
By Theorem 4.4.1 in Kirch~\cite{kirchdiss}, it holds that
%
\begin{equation}\label{eq_pcor_boot_1}\quad
\Biggl|\sum_{j=1}^N\bigl(\cos(t_1\lambda_j)\cos(t_2\lambda_j)+\sin
(t_1\lambda_j)\sin(t_2\lambda_j)\bigr)\Biggr|\ls
\cases{N, &\quad $t_1=t_2$,\cr
1, &\quad $t_1\neq t_2$.}
\end{equation}
By (\ref{cor_boot_est}) and an application of the Cauchy--Schwarz
inequality this implies
%
\begin{eqnarray}\label{eq_pcor_boot_4}
|F_T(j)|&\ls&\sum_{t=1}^T|V(t)-\wv(t)|=o_P( T\alpha_T^{-1/2}
),\nonumber\\[-2pt]
\sum_{j=1}^{2N}F^2_T(j)&\preceq&
N\sum_{t=1}^T\bigl(V(t)-\wv(t)\bigr)^2\nonumber\\[-9pt]\\[-9pt]
&&{}+\sum
_{t_1\neq t_2}\bigl|\bigl(V(t_1)-\wv(t_1)\bigr)\bigl(V(t_2)-\wv(t_2)\bigr)\bigr|\nonumber\\[-2pt]
&=&o_P( T^2\alpha_T^{-1} ).\nonumber
\end{eqnarray}
Equation (\ref{eq_pcor_boot_3}) follows by
%
\begin{eqnarray}\label{eq_pcor_boot_2}
x_V(j)-x_{\wv}(j)
&=& T^{-1/2} F_T(j), \nonumber\\[-9pt]\\[-9pt] y_V(j)-y_{\wv}(j)&=& T^{-1/2}
F_T(N+j).\nonumber
\end{eqnarray}
\upqed\end{pf*}
\begin{pf*}{Proof of Lemma~\ref{lem_spectralestimate}}
Some calculations show that the lemma essentially rephrases Theorem 2.1
in Robinson~\cite{robinson91}, respectively, Theorem 3.2 in Shao and
Wu~\cite{shaowu07}.
\end{pf*}
\begin{pf*}{Proof of Lemma~\ref{lem_llnperio}}
Some careful calculations yield
%
\begin{eqnarray} \label{eq_cov_3}
\sup_{1\ls l,k\ls N}| \cov(x(l),x(k))-\pi f(\lambda_k)\delta
_{l,k}|&\to&0,\nonumber\\[-9pt]\\[-9pt]
\sup_{1\ls l,k\ls N}| \cov(y(l),y(k))-\pi f(\lambda_k)\delta
_{l,k}|&\to&0.\nonumber
\end{eqnarray}
This implies by $\E x(k)=\E y(k)=0$ and an application of the
Markov inequality yields
\[
\frac{1}{2N}\sum_{j=1}^N\frac{x(j)}{\sqrt{f(\lambda
_j)}}=o_P(1),\qquad \frac{1}{2N}\sum_{j=1}^N\frac{y(j)}{\sqrt
{f(\lambda_j)}}=o_P(1),
\]
hence assertion (a). Similar arguments using Proposition 10.3.1 in
Brockwell and Davis~\cite{brockwelldavis} yield assertions (b) and (c).
\end{pf*}
\begin{pf*}{Proof of Lemma~\ref{lem_edf}}
Similar arguments as in the proof of Corollary~2.2 in Shao and Wu \cite
{shaowu07} yield the result. Additionally an argument similar to Freedman and
Lane~\cite{freedmanlane80} is needed.
\end{pf*}
\begin{pf*}{Proof of Lemma~\ref{lem_llnperio2}}
Similar arguments as in the proof of Theorem~A.1 in Franke and H\"
{a}rdle~\cite{franke} yield the result.
\end{pf*}
\begin{pf*}{Proof of Theorem~\ref{th_cpa}}
It is sufficient to prove the assertion of Corollary~\ref
{cor_boot_est} under $H_0$ as well as $H_1$, then the assertion follows
from Theorem~\ref{th_main} as well as\vadjust{\goodbreak} the continuous mapping theorem.
By the H\'{a}jek--Renyi inequality it follows under $H_0$ that 
%
\begin{eqnarray*}
\frac{1}{T}\sum_{t=1}^T\bigl(V(t)-\widehat{V}(t)\bigr)^2&=&\frac{\widehat
{\widetilde{k}}}{T}(\mu-\widehat{\mu}_1)^2+\frac{T-\widehat{\widetilde
{k}}}{T}(\mu-\widehat{\mu}_2)^2\\[-2pt]
&=&\frac{\log T}{T}\Biggl(\frac{1}{\sqrt{(\log T) \widehat{\widetilde
{k}}}} \sum_{j=1}^{\widehat{\widetilde{k}}}\bigl(V(t)-\E(V(t))\bigr) \Biggr)^2\\[-2pt]
&&{} +\frac{\log T}{T}\Biggl(\frac{1}{\sqrt{(\log T) (T-\widehat
{\widetilde{k}})}} \sum_{j=\widehat{\widetilde{k}}+1}^T\bigl(V(t)-\E(V(t))\bigr)
\Biggr)^2\\[-2pt]
&=&O_P\biggl( \frac{\log T}{T} \biggr),
\end{eqnarray*}
which yields the assertion of Corollary~\ref{cor_boot_est}. Similarly,
under alternatives
\begin{eqnarray*}
&&\frac{1}{T}\sum_{t=1}^T\bigl(V(t)-\widehat{V}(t)\bigr)^2\\[-2pt]
&&\qquad=
\frac{\min(\widehat{\widetilde{k}},\widetilde{k})}{T}
(\mu_1-\widehat{\mu}_1)^2+ |d+\mu_j-\widehat{\mu}_j|^2 \frac
{|\widehat{\widetilde{k}}-{\widetilde{k}}|}{T}\\[-2pt]
&&\qquad\quad{} + \frac{T-\max(\widehat
{\widetilde{k}},\widetilde{k})}{T} (\mu_2-\widehat{\mu}_2)^2\\[-2pt]
&&\qquad=O_P\biggl(\max\biggl( \frac{\log T}{T}, \beta_T\biggr) \biggr),
\end{eqnarray*}
where $d=\mu_1-\mu_2$ and $j=2$ if $\widehat{\widetilde{k}}<k$ and
$d=\mu_2-\mu_1$ and $j=1$ otherwise, which yields the assertion of
Corollary~\ref{cor_boot_est}.
\end{pf*}
%
%
\begin{pf*}{Proof of Theorem~\ref{th_unitroot}}
Noting that $Y^*(k)=\sum_{j=1}^k V^*(j)$, the assertion follows from
an application of Corollaries~\ref{cor_boot_est} and~\ref{cor_main}
as well as (\ref{eq_cons_boot_tau}), since $V(t)-\widehat{V}(t)=(\rho
-\widehat{\rho}_T)Y(t-1)$.
\end{pf*}

\begin{supplement}
\stitle{Detailed proofs}
\slink[doi]{10.1214/10-AOS868SUPP} 
\sdatatype{.pdf}
\sfilename{aos868\_suppl.pdf}
\sdescription{In this supplement we give the detailed technical proofs
of the previous sections.}
\end{supplement}


%

%
\printaddresses


\begin{thebibliography}{59}

\bibitem{anthuslp}
\begin{barticle}[author]
\bauthor{\bsnm{Antoch},~\bfnm{J.}\binits{J.}},
  \bauthor{\bsnm{Hu{\v{s}}kov{\'{a}}},~\bfnm{M.}\binits{M.}} \AND
  \bauthor{\bsnm{Pr{\'{a}}{\v{s}}kov{\'{a}}},~\bfnm{Z.}\binits{Z.}}
(\byear{1997}).
\btitle{Effect of dependence on statistics for determination of change}.
\bjournal{J. Statist. Plann. Inference}
\bvolume{60}
\bpages{291--310}.
\end{barticle}
\MR{1456633}
\endbibitem

\bibitem{billv2}
\begin{bbook}[author]
\bauthor{\bsnm{Billingsley},~\bfnm{P.}\binits{P.}}
(\byear{1999}).
\btitle{Convergence of Probability Measures}, \bedition{2nd} ed.
\bpublisher{Wiley}, \baddress{New York}.
\end{bbook}
\MR{1700749}\vadjust{\goodbreak}
\endbibitem

\bibitem{braun}
\begin{barticle}[author]
\bauthor{\bsnm{Braun},~\bfnm{W.~J.}\binits{W.~J.}} \AND
  \bauthor{\bsnm{Kulperger},~\bfnm{R.~J.}\binits{R.~J.}}
(\byear{1997}).
\btitle{Properties of a Fourier bootstrap method for time series}.
\bjournal{Comm. Statist. Theory Methods}
\bvolume{26}
\bpages{1329--1336}.
\end{barticle}
\MR{1456834}
\endbibitem

\bibitem{brillinger}
\begin{bbook}[author]
\bauthor{\bsnm{Brillinger},~\bfnm{D.~R.}\binits{D.~R.}}
(\byear{1981}).
\btitle{Time Series. Data Analysis and Theory},
\bedition{2nd} ed.
\bpublisher{Holden-Day}, \baddress{San Francisco}.
\end{bbook}
\MR{0595684}
\endbibitem

\bibitem{brockwelldavis}
\begin{bbook}[author]
\bauthor{\bsnm{Brockwell},~\bfnm{P.~J.}\binits{P.~J.}} \AND
  \bauthor{\bsnm{Davis},~\bfnm{R.~A.}\binits{R.~A.}}
(\byear{1991}).
\btitle{Time Series: Theory and Methods}, \bedition{2nd} ed.
\bpublisher{Springer}, \baddress{New York}.
\end{bbook}
\MR{1093459}
\endbibitem

\bibitem{buehlmann02}
\begin{barticle}[author]
\bauthor{\bsnm{B{\"{u}}hlmann},~\bfnm{P.}\binits{P.}}
(\byear{2002}).
\btitle{Bootstraps for time series}.
\bjournal{Statist. Sci.}
\bvolume{17}
\bpages{52--72}.
\end{barticle}
\MR{1910074}
\endbibitem

\bibitem{cavataylor08}
\begin{barticle}[author]
\bauthor{\bsnm{Cavaliere},~\bfnm{G.}\binits{G.}} \AND
  \bauthor{\bsnm{Taylor},~\bfnm{A.~M.~R.}\binits{A.~M.~R.}}
(\byear{2008}).
\btitle{Bootstrap unit root tests for time series with non-stationary
  volatility.}
\bjournal{Econom. Theory}
\bvolume{24}
\bpages{43--71}.
\end{barticle}
\MR{2408858}
\endbibitem

\bibitem{cavataylor09a}
\begin{barticle}[author]
\bauthor{\bsnm{Cavaliere},~\bfnm{G.}\binits{G.}} \AND
  \bauthor{\bsnm{Taylor},~\bfnm{A.~M.~R.}\binits{A.~M.~R.}}
(\byear{2009}).
\btitle{Bootstrap $M$ unit root tests.}
\bjournal{Econometric Rev.}
\bvolume{28}
\bpages{393--421}.
\end{barticle}
\MR{2555314}
\endbibitem

\bibitem{cavataylor09b}
\begin{barticle}[author]
\bauthor{\bsnm{Cavaliere},~\bfnm{G.}\binits{G.}} \AND
  \bauthor{\bsnm{Taylor},~\bfnm{A.~M.~R.}\binits{A.~M.~R.}}
(\byear{2009}).
\btitle{Heteroskedastic time series with a unit root}.
\bjournal{Econom. Theory}
\bvolume{25}
\bpages{379--400}.
\end{barticle}
\MR{2540499}
\endbibitem

\bibitem{chan97}
\begin{barticle}[author]
\bauthor{\bsnm{Chan},~\bfnm{K.~S.}\binits{K.~S.}}
(\byear{1997}).
\btitle{On the validity of the method of surrogate data}.
\bjournal{Fields Inst. Commun.}
\bvolume{11}
\bpages{77--97}.
\end{barticle}
\MR{1426615}
\endbibitem

\bibitem{changpark03}
\begin{barticle}[author]
\bauthor{\bsnm{Chang},~\bfnm{Y.}\binits{Y.}} \AND
  \bauthor{\bsnm{Park},~\bfnm{J.~Y.}\binits{J.~Y.}}
(\byear{2003}).
\btitle{A sieve bootstrap for the test of a unit root}.
\bjournal{J.~Time Ser. Anal.}
\bvolume{24}
\bpages{379--400}.
\end{barticle}
\MR{1997120}
\endbibitem

\bibitem{chiu88}
\begin{barticle}[author]
\bauthor{\bsnm{Chiu},~\bfnm{S.~T.}\binits{S.~T.}}
(\byear{1988}).
\btitle{{Weighted least squares estimators on the frequency domain for the
  parameters of a time series}}.
\bjournal{Ann. Statist.}
\bvolume{16}
\bpages{1315--1326}.
\end{barticle}
\MR{0959204}
\endbibitem

\bibitem{csoehor2}
\begin{bbook}[author]
\bauthor{\bsnm{Cs{\"{o}}rg{\H{o}}},~\bfnm{M.}\binits{M.}} \AND
  \bauthor{\bsnm{Horv{\'{a}}th},~\bfnm{L.}\binits{L.}}
(\byear{1993}).
\btitle{Weighted Approximations in Probability and Statistics}.
\bpublisher{Wiley}, \baddress{Chichester}.
\end{bbook}
\MR{1215046}
\endbibitem

\bibitem{csoehor}
\begin{bbook}[author]
\bauthor{\bsnm{Cs{\"{o}}rg{\H{o}}},~\bfnm{M.}\binits{M.}} \AND
  \bauthor{\bsnm{Horv{\'{a}}th},~\bfnm{L.}\binits{L.}}
(\byear{1997}).
\btitle{Limit Theorems in Change-Point Analysis}.
\bpublisher{Wiley}, \baddress{Chichester}.
\end{bbook}
\MR{2743035}
\endbibitem

\bibitem{dahlhausjanas}
\begin{barticle}[author]
\bauthor{\bsnm{Dahlhaus},~\bfnm{R.}\binits{R.}} \AND
  \bauthor{\bsnm{Janas},~\bfnm{D.}\binits{D.}}
(\byear{1996}).
\btitle{A frequency domain bootstrap for ratio statistics in time series
  analysis}.
\bjournal{Ann. Statist.}
\bvolume{24}
\bpages{1934--1963}.
\end{barticle}
\MR{1421155}
\endbibitem

\bibitem{dickeyfuller79}
\begin{barticle}[author]
\bauthor{\bsnm{Dickey},~\bfnm{D.~A.}\binits{D.~A.}} \AND
  \bauthor{\bsnm{Fuller},~\bfnm{W.~A.}\binits{W.~A.}}
(\byear{1979}).
\btitle{Distribution of the estimators for autoregressive time series with a
  unit root}.
\bjournal{J. Amer. Statist. Assoc.}
\bvolume{74}
\bpages{427--431}.
\end{barticle}
\MR{0548036}
\endbibitem

\bibitem{efron79}
\begin{barticle}[author]
\bauthor{\bsnm{Efron},~\bfnm{B.}\binits{B.}}
(\byear{1979}).
\btitle{Bootstrap methods: {A}nother look at the jackknife}.
\bjournal{Ann. Statist.}
\bvolume{7}
\bpages{1--26}.
\end{barticle}
\MR{0515681}
\endbibitem

\bibitem{franke}
\begin{barticle}[author]
\bauthor{\bsnm{Franke},~\bfnm{J.}\binits{J.}} \AND
  \bauthor{\bsnm{H{\"{a}}rdle},~\bfnm{W.}\binits{W.}}
(\byear{1992}).
\btitle{On bootstrapping kernel spectral estimates}.
\bjournal{Ann. Statist.}
\bvolume{20}
\bpages{121--145}.
\end{barticle}
\MR{1150337}
\endbibitem

\bibitem{freedmanlane80}
\begin{barticle}[author]
\bauthor{\bsnm{Freedman},~\bfnm{D.}\binits{D.}} \AND
  \bauthor{\bsnm{Lane},~\bfnm{D.}\binits{D.}}
(\byear{1980}).
\btitle{The empirical distribution of {F}ourier coefficients}.
\bjournal{Ann. Statist.}
\bvolume{8}
\bpages{1244--1251}.
\end{barticle}
\MR{0594641}
\endbibitem

\bibitem{fuller96}
\begin{bbook}[author]
\bauthor{\bsnm{Fuller},~\bfnm{W.~A.}\binits{W.~A.}}
(\byear{1996}).
\btitle{Introduction to Statistical Time Series}, \bedition{2nd} ed.
\bpublisher{Wiley}, \baddress{New York}.
\end{bbook}
\MR{1365746}
\endbibitem

\bibitem{goetzekuensch96}
\begin{barticle}[author]
\bauthor{\bsnm{G{\"{o}}tze},~\bfnm{F.}\binits{F.}} \AND
  \bauthor{\bsnm{K{\"{u}}nsch},~\bfnm{H.~R.}\binits{H.~R.}}
(\byear{1996}).
\btitle{Second-order correctness of the blockwise bootstrap for stationary
  observations}.
\bjournal{Ann. Statist.}
\bvolume{24}
\bpages{1914--1933}.
\end{barticle}
\MR{1421154}
\endbibitem

\bibitem{hamilton94}
\begin{bbook}[author]
\bauthor{\bsnm{Hamilton},~\bfnm{J.}\binits{J.}}
(\byear{1994}).
\btitle{Time Series Analysis}.
\bpublisher{Princeton Univ. Press}, \baddress{Princeton, NJ}.
\end{bbook}
\MR{1278033}
\endbibitem

\bibitem{hidalgo}
\begin{barticle}[author]
\bauthor{\bsnm{Hidalgo},~\bfnm{J.}\binits{J.}}
(\byear{2003}).
\btitle{An alternative bootstrap to moving blocks for time series regression
  models}.
\bjournal{J. Econometrics}
\bvolume{117}
\bpages{369--399}.
\end{barticle}
\MR{2008775}
\endbibitem

\bibitem{horlp}
\begin{barticle}[author]
\bauthor{\bsnm{Horv{\'{a}}th},~\bfnm{L.}\binits{L.}}
(\byear{1997}).
\btitle{Detection of changes in linear sequences}.
\bjournal{Ann. Inst. Statist. Math.}
\bvolume{49}
\bpages{271--283}.
\end{barticle}
\MR{1463306}
\endbibitem

\bibitem{hurvichzeger}
\begin{bmisc}[author]
\bauthor{\bsnm{Hurvich},~\bfnm{C.~M.}\binits{C.~M.}} \AND
  \bauthor{\bsnm{Zeger},~\bfnm{S.~L.}\binits{S.~L.}}
(\byear{1987}).
\bhowpublished{Frequency domain bootstrap methods for time series.
Working paper, New York Univ.}
\end{bmisc}
\endbibitem

\bibitem{huskirchconfstud}
\begin{barticle}[author]
\bauthor{\bsnm{Hu{\v{s}}kov{\'{a}}},~\bfnm{M.}\binits{M.}} \AND
  \bauthor{\bsnm{Kirch},~\bfnm{C.}\binits{C.}}
(\byear{2010}).
\btitle{A note on studentized confidence intervals in change-point analysis}.
\bjournal{Comput. Statist.}
\bvolume{25}
\bpages{269--289}.
\end{barticle}
\endbibitem

\bibitem{dahlhausjanas2}
\begin{binproceedings}[author]
\bauthor{\bsnm{Janas},~\bfnm{D.}\binits{D.}} \AND
  \bauthor{\bsnm{Dahlhaus},~\bfnm{R.}\binits{R.}}
(\byear{1994}).
\btitle{A frequency domain bootstrap for time series}.
In \bbooktitle{Proceedings of the 26th Symposium on the Interface}
(\beditor{J. Sall and A. Lehman}, eds.)
\bpages{423--425}.
\bpublisher{Interface Foundation of North America},
\baddress{Fairfax Station, VA}.
\end{binproceedings}
\endbibitem

\bibitem{jentschkreiss10}
\begin{barticle}[author]
\bauthor{\bsnm{Jentsch},~\bfnm{C.}\binits{C.}} \AND
  \bauthor{\bsnm{Kreiss},~\bfnm{J.~P.}\binits{J.~P.}}
(\byear{2010}).
\btitle{The multiple hybrid bootstrap---resampling multivariate linear
  processes.}
\bjournal{J. Multivariate Anal.}
\bvolume{101}
\bpages{2320--2345}.
\end{barticle}
\endbibitem

\bibitem{kirchdiss}
\begin{bmisc}[author]
\bauthor{\bsnm{Kirch},~\bfnm{C.}\binits{C.}}
(\byear{2006}).
\bhowpublished{Resampling methods for the change analysis of dependent
data. Ph.D. thesis, Univ. Cologne. Available at
\texttt{\href{http://kups.ub.uni-koeln.de/volltexte/2006/1795/}{http://kups.ub.uni-koeln.de/}
\href{http://kups.ub.uni-koeln.de/volltexte/2006/1795/}{volltexte/2006/1795/}}.}
\end{bmisc}
\endbibitem

\bibitem{kirchblock}
\begin{barticle}[author]
\bauthor{\bsnm{Kirch},~\bfnm{C.}\binits{C.}}
(\byear{2007}).
\btitle{Block permutation principles for the change analysis of dependent
  data}.
\bjournal{J.~Statist. Plann. Inference}
\bvolume{137}
\bpages{2453--2474}.
\end{barticle}
\MR{2325449}
\endbibitem

\bibitem{kirchfreq}
\begin{barticle}[author]
\bauthor{\bsnm{Kirch},~\bfnm{C.}\binits{C.}}
(\byear{2007}).
\btitle{Resampling in the frequency domain of time series to determine critical
  values for change-point tests}.
\bjournal{Statist. Decisions}
\bvolume{25}
\bpages{237--261}.
\end{barticle}
\MR{2412072}
\endbibitem

\bibitem{suppA}
\begin{bmisc}[author]
\bauthor{\bsnm{Kirch},~\bfnm{C.}\binits{C.}} \AND
  \bauthor{\bsnm{Politis},~\bfnm{D.~N.}\binits{D.~N.}}
(\byear{2010}).
\bhowpublished{Supplement to ``{TFT}-{b}ootstrap: {R}esampling time series in the
frequency domain to obtain replicates in the time domain.''
\href{http://dx.doi.org/10.1214/10-AOS868SUPP}{DOI:10.1214/10-AOS868SUPP}.}
\end{bmisc}
\endbibitem

\bibitem{kokleip98}
\begin{barticle}[author]
\bauthor{\bsnm{Kokoszka},~\bfnm{P.}\binits{P.}} \AND
  \bauthor{\bsnm{Leipus},~\bfnm{R.}\binits{R.}}
(\byear{1998}).
\btitle{Change-point in the mean of dependent observations}.
\bjournal{Statist. Probab. Lett.}
\bvolume{40}
\bpages{385--393}.
\end{barticle}
\MR{1664564}
\endbibitem

\bibitem{papakreiss03}
\begin{barticle}[author]
\bauthor{\bsnm{Kreiss},~\bfnm{J.~P.}\binits{J.~P.}} \AND
  \bauthor{\bsnm{Paparoditis},~\bfnm{E.}\binits{E.}}
(\byear{2003}).
\btitle{Autoregressive-aided periodogram bootstrap for time series.}
\bjournal{Ann. Statist.}
\bvolume{31}
\bpages{1923--1955}.
\end{barticle}
\MR{2036395}
\endbibitem

\bibitem{lahiri03a}
\begin{barticle}[author]
\bauthor{\bsnm{Lahiri},~\bfnm{S.~N.}\binits{S.~N.}}
(\byear{2003}).
\btitle{A necessary and sufficient condition for asymptotic independence of
  discrete Fourier transforms under short- and long-range dependence}.
\bjournal{Ann. Statist.}
\bvolume{31}
\bpages{613--641}.
\end{barticle}
\MR{1983544}
\endbibitem

\bibitem{lahiri03}
\begin{bbook}[author]
\bauthor{\bsnm{Lahiri},~\bfnm{S.~N.}\binits{S.~N.}}
(\byear{2003}).
\btitle{Resampling Methods for Dependent Data}.
\bpublisher{Springer}, \baddress{New York}.
\end{bbook}
\MR{2001447}
\endbibitem

\bibitem{mammen08}
\begin{binproceedings}[author]
\bauthor{\bsnm{Maiwald},~\bfnm{T.}\binits{T.}},
  \bauthor{\bsnm{Mammen},~\bfnm{E.}\binits{E.}},
  \bauthor{\bsnm{Nandi},~\bfnm{S.}\binits{S.}} \AND
  \bauthor{\bsnm{Timmer},~\bfnm{J.}\binits{J.}}
(\byear{2008}).
\btitle{Surrogate data---A~qualitative and quantitative analysis}.
In \bbooktitle{Mathematical Methods in Time Series Analysis and Digital Image
  Processing}
(\beditor{R. Dahlhaus, J. Kurths, P. Maas and
J. Timmer}, eds.).
\bpublisher{Springer}, \baddress{Berlin}.
\end{binproceedings}
\endbibitem

\bibitem{mallows72}
\begin{barticle}[author]
\bauthor{\bsnm{Mallows},~\bfnm{C.~L.}\binits{C.~L.}}
(\byear{1972}).
\btitle{A note on asymptotic joint normality}.
\bjournal{Ann. Math. Statist.}
\bvolume{43}
\bpages{508--515}.
\end{barticle}
\MR{0298812}
\endbibitem

\bibitem{mammennandipp}
\begin{barticle}[author]
\bauthor{\bsnm{Mammen},~\bfnm{E.}\binits{E.}} \AND
  \bauthor{\bsnm{Nandi},~\bfnm{S.}\binits{S.}}
(\byear{2008}).
\btitle{Some theoretical properties of phase randomized multivariate
  surrogates}.
\bjournal{Statistics}
\bvolume{42}
\bpages{195--205}.
\end{barticle}
\endbibitem

\bibitem{ngperron01}
\begin{barticle}[author]
\bauthor{\bsnm{Ng},~\bfnm{S.}\binits{S.}} \AND
  \bauthor{\bsnm{Perron},~\bfnm{P.}\binits{P.}}
(\byear{2001}).
\btitle{Lag length selection and the construction of unit root tests with good
  size and power}.
\bjournal{Econometrica}
\bvolume{69}
\bpages{1519--1554}.
\end{barticle}
\MR{1865220}
\endbibitem

\bibitem{papa02}
\begin{binproceedings}[author]
\bauthor{\bsnm{Paparoditis},~\bfnm{E.}\binits{E.}}
(\byear{2002}).
\btitle{Frequency domain bootstrap for time series}.
In \bbooktitle{Empirical Process Techniques for Dependent Data}
(\beditor{H. Dehling et al.}, eds.)
\bpages{365--381}.
\bpublisher{Birkh\"{a}user}, \baddress{Boston, MA}.
\end{binproceedings}
\MR{1958790}
\endbibitem

\bibitem{papapolitis99}
\begin{barticle}[author]
\bauthor{\bsnm{Paparoditis},~\bfnm{E.}\binits{E.}} \AND
  \bauthor{\bsnm{Politis},~\bfnm{D.~N.}\binits{D.~N.}}
(\byear{1999}).
\btitle{The local bootstrap for periodogram statistics}.
\bjournal{J. Time Ser. Anal.}
\bvolume{20}
\bpages{193--222}.
\end{barticle}
\MR{1701054}
\endbibitem

\bibitem{papapolitis03}
\begin{barticle}[author]
\bauthor{\bsnm{Paparoditis},~\bfnm{E.}\binits{E.}} \AND
  \bauthor{\bsnm{Politis},~\bfnm{D.~N.}\binits{D.~N.}}
(\byear{2003}).
\btitle{{Residual-based block bootstrap for unit root testing}}.
\bjournal{Econometrica}
\bvolume{71}
\bpages{813--855}.
\end{barticle}
\MR{1983228}
\endbibitem

\bibitem{papapolitis05}
\begin{barticle}[author]
\bauthor{\bsnm{Paparoditis},~\bfnm{E.}\binits{E.}} \AND
  \bauthor{\bsnm{Politis},~\bfnm{D.~N.}\binits{D.~N.}}
(\byear{2005}).
\btitle{Bootstrapping unit root tests for autoregressive time series}.
\bjournal{J. Amer. Statist. Assoc.}
\bvolume{100}
\bpages{545--553}.
\end{barticle}
\MR{2160558}
\endbibitem

\bibitem{park03}
\begin{barticle}[author]
\bauthor{\bsnm{Park},~\bfnm{J.~Y.}\binits{J.~Y.}}
(\byear{2003}).
\btitle{Bootstrap unit root tests}.
\bjournal{Econometrica}
\bvolume{71}
\bpages{1845--1895}.
\end{barticle}
\MR{2015421}
\endbibitem

\bibitem{perronng96}
\begin{barticle}[author]
\bauthor{\bsnm{Perron},~\bfnm{P.}\binits{P.}} \AND
  \bauthor{\bsnm{Ng},~\bfnm{S.}\binits{S.}}
(\byear{1996}).
\btitle{Useful modifications to some unit root tests with dependent errors and
  their local asymptotic properties}.
\bjournal{Rev. Econom. Stud.}
\bvolume{63}
\bpages{435--463}.
\end{barticle}
\endbibitem

\bibitem{phillips87}
\begin{barticle}[author]
\bauthor{\bsnm{Phillips},~\bfnm{P.~C.~B.}\binits{P.~C.~B.}}
(\byear{1987}).
\btitle{{Time series regression with a unit root}}.
\bjournal{Econometrica}
\bvolume{55}
\bpages{277--301}.
\end{barticle}
\MR{0882096}
\endbibitem

\bibitem{phillipsperron88}
\begin{barticle}[author]
\bauthor{\bsnm{Phillips},~\bfnm{P.~C.~B.}\binits{P.~C.~B.}} \AND
  \bauthor{\bsnm{Perron},~\bfnm{P.}\binits{P.}}
(\byear{1988}).
\btitle{Testing for a unit root in time series regression}.
\bjournal{Biometrika}
\bvolume{75}
\bpages{335--346}.
\end{barticle}
\MR{0946054}
\endbibitem

\bibitem{politis03}
\begin{barticle}[author]
\bauthor{\bsnm{Politis},~\bfnm{D.~N.}\binits{D.~N.}}
(\byear{2003}).
\btitle{Adaptive bandwidth choice}.
\bjournal{J. Nonparametr. Stat.}
\bvolume{15}
\bpages{517--533}.
\end{barticle}
\MR{2017485}
\endbibitem

\bibitem{politis03a}
\begin{barticle}[author]
\bauthor{\bsnm{Politis},~\bfnm{D.~N.}\binits{D.~N.}}
(\byear{2003}).
\btitle{The impact of bootstrap methods on time series analysis}.
\bjournal{Statist. Sci.}
\bvolume{18}
\bpages{219--230}.
\end{barticle}
\MR{2026081}
\endbibitem

\bibitem{politis05}
\begin{bmisc}[author]
\bauthor{\bsnm{Politis},~\bfnm{D.~N.}\binits{D.~N.}}
(\byear{2009}).
\bhowpublished{Higher-order accurate, positive semi-definite estimation of
  large-sample covariance and spectral density matrices.
Preprint. Paper 2005-03R, Dept. Economics, UCSD.
Available at \texttt{\href{http://repositories.cdlib.org/ucsdecon/2005-03R}{http://repositories.cdlib.org/}
\href{http://repositories.cdlib.org/ucsdecon/2005-03R}{ucsdecon/2005-03R}}.}
\end{bmisc}
\endbibitem

\bibitem{priestleybook}
\begin{bbook}[author]
\bauthor{\bsnm{Priestley},~\bfnm{M.~B.}\binits{M.~B.}}
(\byear{1981}).
\btitle{Spectral Analysis and Time Series}.
\bpublisher{Academic Press}, \baddress{London}.
\end{bbook}
\endbibitem

\bibitem{robinson91}
\begin{barticle}[author]
\bauthor{\bsnm{Robinson},~\bfnm{P.~M.}\binits{P.~M.}}
(\byear{1991}).
\btitle{Automatic frequency domain inference on semiparametric and
  nonparametric models}.
\bjournal{Econometrica}
\bvolume{59}
\bpages{1329--1363}.
\end{barticle}
\MR{1133037}
\endbibitem

\bibitem{romanothombs96}
\begin{barticle}[author]
\bauthor{\bsnm{Romano},~\bfnm{J.~P.}\binits{J.~P.}} \AND
  \bauthor{\bsnm{Thombs},~\bfnm{L.~A.}\binits{L.~A.}}
(\byear{1996}).
\btitle{Inference for autocorrelations under weak assumptions.}
\bjournal{J. Amer. Statist. Assoc.}
\bvolume{91}
\bpages{590--600}.
\end{barticle}
\MR{1395728}
\endbibitem

\bibitem{shaowu07}
\begin{barticle}[author]
\bauthor{\bsnm{Shao},~\bfnm{X.}\binits{X.}} \AND
  \bauthor{\bsnm{Wu},~\bfnm{W.~B.}\binits{W.~B.}}
(\byear{2007}).
\btitle{Asymptotic spectral theory for nonlinear time series}.
\bjournal{Ann. Statist.}
\bvolume{35}
\bpages{1773--1801}.
\end{barticle}
\MR{2351105}
\endbibitem

\bibitem{stock99}
\begin{bincollection}[author]
\bauthor{\bsnm{Stock},~\bfnm{J.~H.}\binits{J.~H.}}
(\byear{1999}).
\btitle{A class of tests for integration and cointegration}.
In \bbooktitle{Cointegration, Causality, and Forecasting---Festschrift in
  Honour of Clive W. J. Granger}
(\beditor{R. F. Engle and
H. White}, eds.)
\bpages{137--167}.
\bpublisher{Oxford Univ. Press}, \baddress{New York}.
\end{bincollection}
\endbibitem

\bibitem{theileretal92}
\begin{barticle}[author]
\bauthor{\bsnm{Theiler},~\bfnm{J.}\binits{J.}},
  \bauthor{\bsnm{Eubank},~\bfnm{S.}\binits{S.}},
  \bauthor{\bsnm{Longtin},~\bfnm{A.}\binits{A.}},
  \bauthor{\bsnm{Galdrikan},~\bfnm{B.}\binits{B.}} \AND
  \bauthor{\bsnm{Farmer},~\bfnm{J.~D.}\binits{J.~D.}}
(\byear{1992}).
\btitle{Testing for nonlinearity in time series: The method of surrogate data}.
\bjournal{Phys. D}
\bvolume{58}
\bpages{77--94}.
\end{barticle}
\endbibitem

\bibitem{wu08}
\begin{bmisc}[author]
\bauthor{\bsnm{Wu},~\bfnm{W.~B.}\binits{W.~B.}}
(\byear{2008}).
\bhowpublished{Personal communication}.
\end{bmisc}
\endbibitem

\end{thebibliography}
\end{document}